\documentclass[journal]{IEEEtran}

\newcommand{\beq}{\begin{equation}}
\newcommand{\eeq}{\end{equation}}
\newcommand{\bsq}{\begin{subequations}}
\newcommand{\esq}{\end{subequations}}
\newcommand{\bq}{\begin{eqnarray}}
\newcommand{\eq}{\end{eqnarray}}
\newcommand{\bqn}{\begin{eqnarray*}}
\newcommand{\eqn}{\end{eqnarray*}}
\usepackage{bbm}
\usepackage[pdftex]{graphicx}
\usepackage{enumerate}
\usepackage{enumitem} 

\usepackage{mathptmx}       
\DeclareMathAlphabet{\mathcal}{OMS}{cmsy}{m}{n}
\usepackage{helvet}         
\usepackage{courier}        
\usepackage{type1cm}        
\usepackage{makeidx}         
\usepackage{graphicx}        
\usepackage{url} 
\usepackage{multicol}        
\usepackage[bottom]{footmisc}
\usepackage{ifthen}
\usepackage{subfigure}
\usepackage[usenames,dvipsnames]{color}

\makeindex             

\usepackage{graphics} 
\usepackage{graphicx} 
\usepackage{epsfig} 
\usepackage{epstopdf}
\usepackage{mathptmx} 
\usepackage{times} 
\usepackage[cmex10]{amsmath}

\usepackage{mathtools}  
\usepackage{amssymb}  
\usepackage{amsthm}

\usepackage{amsfonts}
\usepackage{cite}
\usepackage{verbatim}
\usepackage{comment}
\usepackage{algorithmic}
\usepackage{multirow}
\usepackage{cite}
\usepackage{stfloats}
\usepackage{supertabular}
\usepackage{longtable}

\DeclareGraphicsExtensions{.pdf,.png,.jpg,.eps,.ps}
\renewcommand{\arraystretch}{1.2}
\usepackage{arydshln}
\usepackage{multicol}
\usepackage{colortbl}
\theoremstyle{definition}
\newtheorem{theorem}{Theorem}

\newtheorem{proposition}{Proposition}

\theoremstyle{definition}

\ifCLASSINFOpdf
\else
\fi

\usepackage{comment}
\usepackage{ifthen}
\newboolean{showcomments}
\setboolean{showcomments}{true}
\newcommand{\ychen}[1]{\ifthenelse{\boolean{showcomments}}
        { \textcolor{red}{YC: #1}}}
\newcommand{\tongxin}[1]{\ifthenelse{\boolean{showcomments}}
        { \textcolor{blue}{(#1)}}{}}

\usepackage{amsmath}
\allowdisplaybreaks[4]

\usepackage[ruled]{algorithm2e}

\begin{document}

%
\title{Robust Generation Dispatch with Strategic Renewable Power Curtailment and Decision-Dependent Uncertainty}
%
%
%

\author{Yue Chen and Wei Wei
\thanks{Y. Chen is with the Department of Mechanical and Automation Engineering, the Chinese University of Hong Kong, Hong Kong SAR (email: yuechen@mae.cuhk.edu.hk)}
\thanks{W. Wei is with the the State Key Laboratory of Power Systems, Department of Electrical Engineering, Tsinghua University, Beijing 100084 China. (wei-wei04@mails.tsinghua.edu.cn)}
\thanks{This work has been submitted to the IEEE for possible publication.
Copyright may be transferred without notice, after which this version may no
longer be accessible.}}
%
%

\markboth{Journal of \LaTeX\ Class Files,~Vol.~XX, No.~X, Feb.~2019}%
{Shell \MakeLowercase{\textit{et al.}}: Bare Demo of IEEEtran.cls for IEEE Journals}
%



\maketitle

\begin{abstract}
As renewable energy sources replace traditional power sources (such as thermal generators), uncertainty grows while there are fewer controllable units. To reduce operational risks and avoid frequent real-time emergency controls, a preparatory schedule of renewable generation curtailment is required. This paper proposes a novel two-stage robust generation dispatch (RGD) model, where the preparatory curtailment schedule is optimized in the pre-dispatch stage. The curtailment schedule will then influence the variation range of real-time renewable power output, resulting in a decision-dependent uncertainty (DDU) set. In the re-dispatch stage, the controllable units adjust their outputs within the reserve capacities to maintain power balancing. To overcome the difficulty in solving the RGD with DDU, an adaptive column-and-constraint generation (AC\&CG) algorithm is developed. We prove that the proposed algorithm can generate the optimal solution in finite iterations. Numerical examples show the advantages of the proposed model and algorithm, and validate their practicability and scalability.
\end{abstract}


%
\IEEEpeerreviewmaketitle

\section*{Nomenclature}
\addcontentsline{toc}{section}{Nomenclature}
\subsection{Abbreviations}
\begin{IEEEdescription}[\IEEEusemathlabelsep\IEEEsetlabelwidth{${\underline P _{mn}}$,${\overline P _{mn}}$}]
\item[C\&CG] Column-and-constraint generation.
\item[DDU] Decision-dependent uncertainty.
\item[DRO] Distributionally robust optimization.
\item[RG] Renewable generator.
\item[RGD] Robust generation dispatch.
\item[RO] Robust optimization.
\item[SO] Stochastic optimization.
\end{IEEEdescription}

\subsection{Indices, Sets, and Functions}
\begin{IEEEdescription}[\IEEEusemathlabelsep\IEEEsetlabelwidth{${\underline P _{mn}}$,${\overline P _{mn}}$}]
\item[$\mathcal{AC}_n$] Set of active constraints in iteration $n$.
\item[$Q(p,r^{\pm},w)$] The optimal objective value of the re-dispatch problem when the RG output is $w$.
\item[$S(x,\xi)$] The optimal objective value of the inner ``max-min'' problem given the first-stage decisions $x,\xi$.
\item[$S_f(x,\xi)$] The optimal objective value of the feasibility-check problem given the first-stage $x,\xi$.
\item[$\mathcal{W}(\xi)$]   Decision-dependent uncertainty set.
\end{IEEEdescription}

\subsection{Parameters}
\begin{IEEEdescription}[\IEEEusemathlabelsep\IEEEsetlabelwidth{${\underline P _{mn}}$,${\overline P _{mn}}$}]
\item[$d_{i}^{\pm}$] Upward/Downward regulation cost coefficient of thermal unit $i$.
\item[$\mathcal{D}_{lt}$] Demand of load $l$ in period $t$.
\item[$F_k$] Capacity of transmission line $k$.
\item[$N_g$] Number of thermal units.
\item[$N_k$] Number of transmission lines.
\item[$N_l$] Number of loads.
\item[$N_w$] Number of renewable generators.
\item[$P_i^u, P_i^l$] Maximum/Minimum output of  thermal unit $i$.
\item[$R_i^{\pm}$] Maximum/Minimum reserve of thermal unit $i$.
\item[$\mathcal{R}_i^d, \mathcal{R}_i^u$] Ramping limit of thermal unit $i$.
\item[$T$] Number of dispatch periods.
\item[$u_{it}$] Binary variable, $u_{it}=1$ if the thermal unit $i$ is on in period $t$; otherwise $u_{it}=0$.
\item[$W_{jt}^e$] Forecast output of RG $j$ in period $t$.
\item[$W_{jt}^u, W_{jt}^l$] Upper/Lower bound of the confidence interval of power output of RG $j$ in period $t$.
\item[$\alpha_i$] Running cost coefficient of thermal unit $i$.
\item[$\beta_i^{\pm}$] Upward/Downward reserve cost coefficient of thermal unit $i$.
\item[$\gamma_j$] Predispatch curtailment cost coefficient of RG $j$.
\item[$\hat \gamma_j$] Redispatch curtailment cost coefficient of RG $j$.
\item[$\pi_{ik},\pi_{jk},\pi_{lk}$] Line flow distribution factors.
\item[$\Gamma_S, \Gamma_T$] Spatial/Temporal uncertainty budget.
\end{IEEEdescription}

\subsection{Decision Variables}
\begin{IEEEdescription}[\IEEEusemathlabelsep\IEEEsetlabelwidth{${\underline P _{mn}}$,${\overline P _{mn}}$}]
\item[$p_{it}$] Contemporary output of thermal unit $i$ in $t$.
\item[$p_{it}^{\pm}$] Incremental outputs of thermal unit $i$ in $t$.
\item[$r_{it}^{\pm}$] Reserve capacity of thermal unit $i$ in period $t$.
\item[$\hat w_{jt}$] Output before pre-dispatch curtailment of RG $j$ in period $t$.
\item[$w_{jt}$] Ouput after pre-dispatch curtailment of RG $j$ in period $t$.
\item[$w_{jt}^u, w_{jt}^l$] The amount by which the power output of RG $j$ deviates upward/downward from the forecast.
\item[$z_{jt},\hat z_{jt}$] Binary variables for linearization.
\item[$\xi_{jt}$] Pre-dispatch curtailment of RG $j$ in period $t$.
\item[$\hat \xi_{jt}$] Re-dispatch curtailment of RG $j$ in period $t$.
\end{IEEEdescription}

\section{Introduction}
%
%
%
%

\IEEEPARstart{O}{ne} effective way to mitigate global warming is the massive deployment of renewable energy. In 2020, wind power with an incremental capacity of 93 GW was installed, boosting the total global capacity to 743 GW \cite{council2021gwec}. Despite a decrease in energy consumption due to COVID-19 in 2020, the global cumulative solar capacity increased by 22\%, establishing a new milestone for the solar sector \cite{schmela2021global}. At the same time, increasing renewable energy penetration poses a significant challenge to power system operations due to its intermittent, volatile, and uncertain characteristics. Generation dispatch (a generalization of economic dispatch) is widely used to deal with uncertainties by altering the real-time output of controllable generator within reserve capacity to hedge against renewable power output fluctuations \cite{wei2014robust}. 

The majority of existing research studies the generation dispatch problem under uncertainty using stochastic optimization (SO) \cite{fu2015multiobjective,hu2020bayesian}, robust optimization (RO) \cite{wei2014robust,lorca2014adaptive}, and distributionally robust optimization (DRO) \cite{chen2016distributionally,wei2015distributionally}. Stochastic optimization aims to minimize the expected cost or to meet certain chance constraints based on the assumption that the renewable generator (RG) output follows a known probability distribution. However, in practice, though we can obtain a rough histogram via historical data, it is hard to procure the exact probability distribution. Deviation from the real distribution may result in sub-optimal or even infeasible solutions. In contrast, RO considers all possible realizations of RG output in a prespecified uncertainty set and comes up with the dispatch strategy that minimizes the cost under the worst-case scenario. Various effective algorithms, such as Bender's 
decomposition \cite{bertsimas2012adaptive} and column-and-constraint generation (C\&CG) \cite{zeng2013solving}, were put forward to get the robust optimal solution efficiently. DRO is in-between SO and RO. Instead of using an uncertainty set, it describes the RG output by all possible probability distributions with the same mean and variance recovered from the historical data. Though DRO is less conservative than RO, it is time-consuming since it involves solving a semi-definite program \cite{wei2015distributionally}. Therefore, this paper adopts the robust generation dispatch (RGD) model.

In most existing literature, the uncertainty set in an RGD model is determined based on all available renewable generation, implicitly assuming that 100\% of renewable generation can be used without curtailment. However, as renewable energy sources replace conventional power sources (e.g., thermal generators), uncertainty increases while there are fewer controllable units. As a result, curtailment of renewable energy is inevitable to maintain system security \cite{liu2018curtailment}. References \cite{wu2016solution} and \cite{morales2018robust} developed two-stage unit commitment models taking into account the wind power curtailment. However, the curtailment only happens in real-time without preparatory schedules. This paper focuses on future power systems with a high penetration of renewable energy. Real-time curtailment only could lead to the immediate shutdown of a large number of renewable generators, which may worsen system dynamic performance or trigger instability/oscillation \cite{yang2021chance}. Moreover, with more renewable generators and fewer thermal units, the power system inertia decreases, which will aggravate the risk of instability \cite{abo2013impacts}. To deal with this challenge, in Spain, both the scheduled curtailment that is determined before the day-ahead market closes and the real-time curtailment have been applied \cite{rogers2010examples}. Reference \cite{yang2021chance} also emphasized the need for considering curtailment in the pre-dispatch stage so that we can have a relatively long time to make adjustments. Furthermore, the preparatory schedule of curtailment can help lower the amount of reserve required since the system is less unpredictable. For the reasons above, in this paper, we take into account preparatory curtailment in the pre-dispatch stage.


\begin{table*}[t]
        \renewcommand{\arraystretch}{1.3}
        \renewcommand{\tabcolsep}{0.9em}
        \centering
        \caption{Summary of recent work on DDU}
        \label{tab:summary}
        \begin{tabular}{ccccc}
                \hline 
              \multicolumn{2}{c}{Model}  & Reference  & Types of DDU & Solution Methodology\\
                \hline
              \multicolumn{2}{c}{\multirow{3}{3cm}{Stochastic optimization}} & \cite{yin2019multi} & Type-1 & Non-anticipativity constraints \\
              & & \cite{tarhan2009stochastic} & Type-1 & Non-anticipativity constraints \\
              & & \cite{zhan2016generation} & Type-2 & A quasi-exact solution approach\\
              & & & & \\
              \multirow{7}{3cm}{Robust optimization} & \multirow{2}{2cm}{Static robust} & \cite{nohadani2018optimization} & Type-2 & Equivalent reformulation \\
              & & \cite{lappas2018robust} & Type-2 & Equivalent reformulation \\
              & & & &\\
              & \multirow{5}{2cm}{Adaptive robust} & \cite{vayanos2020robust} & Type-1 & K-adaptability approximation \\
              & & \cite{avraamidou2020adjustable} & Type-2 & Multi-parametric programming \\
              & & \cite{zhang2021robust} & Type-2 & Modified Bender's decomposition \\
              & & \cite{zhang2022two} & Type-2 & Modified Bender's decomposition \\
              & & \cite{su2020robust} & Type-2 & Improved C\&CG\\
              \hline
        \end{tabular}
\end{table*}

When the curtailment of renewable energy is incorporated in the pre-dispatch stage, the variation range of real-time RG output will be influenced. Therefore, distinct from the uncertainty sets in traditional robust optimization models that are pre-specified, the uncertainty set in our model depends on the first-stage decision. \emph{Decision-dependent uncertainty} (DDU) has captured great attention in recent years and mainly two types of DDU have been considered \cite{hellemo2018decision}: 1) information structure related DDU (Type-1), i.e., decisions affect when knowledge about uncertain factors becomes available; 2) uncertainty distribution related DDU (Type-2), i.e., decisions affect the size of an uncertainty set or the chance that a scenario happens. Related studies are summarized in TABLE \ref{tab:summary}. Stochastic programming was applied to deal with DDU \cite{yin2019multi,tarhan2009stochastic,zhan2016generation}. As for robust optimization, reference \cite{nohadani2018optimization} provided a reformulation method to solve the static robust model with DDU in which only ``here-and-now'' variables were included and the size of uncertainty set merely depends on binary variables. Another robust counterpart model was established in \cite{lappas2018robust} with a more general set. In practice, adjustable robust optimization is more commonly used since it can describe the adjustment behaviors after the realization of uncertain factors. The solution of an adjustable robust optimization with DDU can be estimated via a K-adaptability approximation based algorithm \cite{vayanos2020robust}. Multi-parametric programming \cite{avraamidou2020adjustable}, modified Bender’s decomposition \cite{zhang2021robust,zhang2022two}, and improved C\&CG \cite{su2020robust} were deployed to offer an exact solution. However, they are either time-consuming with a growing number of units or concentrated on the linear DDU set. In general, robust optimization with DDU is an important albeit challenging problem due to its computational intractability \cite{zeng2022two}.


This paper proposes a robust generation dispatch model considering curtailment in the pre-dispatch stage. An adaptive C\&CG algorithm is developed to solve the problem efficiently in finite iterations. Our main contributions are two-fold:

1) \textbf{Robust generation dispatch model with DDU}. A novel RGD model considering renewable power curtailment in the pre-dispatch stage is proposed. The pre-dispatch curtailment is defined as the upper limit of the actual RG output, which is co-optimized with the contemporary output and reserve capacity of thermal generators. This enables a preparatory schedule of curtailment and avoids frequent real-time emergency controls, which is more practical. Then in the re-dispatch stage, if the reserve is not enough to maintain power balancing, real-time curtailment will be executed. Distinct from the traditional RGD models, the variation range of RG output is capped by the pre-dispatch curtailment strategy, and so the uncertainty set is decision-dependent. Therefore, the proposed model renders a robust optimization with DDU. The nonlinear DDU set is further linearized into a mixed-integer linear set. Case studies show that the proposed model can avoid significant wind curtailment in the re-dispatch stage and can also lower the amount of reserve required.

2) \textbf{Adaptive C\&CG Algorithm}. With decision-dependent uncertainty, the traditional robust optimization algorithms cannot be applied since the previously selected scenario may be outside of the uncertainty set as the pre-dispatch strategy changes. In this paper, an adaptive C\&CG algorithm is developed. We add the active constraints under the worst-case scenario rather than the worst-case scenario itself to the master problem in each iteration. Two types of degeneracy that may occur during active constraints identification are addressed. When the pre-dispatch strategy changes, the scenario generated by the active constraints remains a vertex of the new uncertainty set. We prove that our algorithm can generate the robust optimal solution in finite iterations. Case studies show that our algorithm can greatly reduce the computational time compared with the nested C\&CG (which may still be applied but without a theoretical guarantee of finite-step-ending).

The rest of this paper is organized as follows. Section \ref{sec-2} introduces the decision-dependent uncertainty set and the robust generation dispatch model. An adaptive C\&CG algorithm is developed in Section \ref{sec-3} to solve the problem. 
Case studies are given in Section \ref{sec-4} with conclusions in Section \ref{sec-5}.

\section{Mathematical Formulation}
\label{sec-2}
We consider a two-stage generation dispatch problem for a system with $N_g$ thermal generators, $N_w$ renewable generators (RGs), and $T$ dispatch periods. In the first-stage (\emph{pre-dispatch stage}), the contemporary output and reserve capacity of thermal generators are determined based on the RG forecast outputs. Moreover, the preparatory curtailment of RGs is also considered to adapt to the increasing fluctuation due to higher renewable energy penetration. In the second-stage (\emph{re-dispatch stage}), the operator adjusts the output of thermal generators within their reserve capacities and/or make additional real-time curtailments to maintain power balancing. In the following, we first introduce the RG output uncertainty set, which is influenced by the first-stage curtailment strategy. Then the robust generation dispatch model with DDU is proposed.
\subsection{Decision-dependent uncertainty set}
In this paper, preparatory curtailment of RGs is considered in the first-stage to avoid frequent calls for emergency controls in real-time. The curtailment is defined as an upper limit of the actual RG output, denoted as $\xi_{jt},\forall j,\forall t$. Denote the actual output and output after preparatory curtailment of RG $j$ in period $t$ as $\hat w_{jt}$ and $w_{jt}$, respectively. Then, the uncertainty set $\mathcal{W}(\xi)$ is as follows, depending on the curtailment strategy $\xi$.
\begin{align}
\label{eq:DDU-set}
    \mathcal{W}(\xi)=\left\{ {w ~\left| {\begin{array}{*{20}{c}}
    w_{jt} = \mbox{min}\{\hat w_{jt} , \xi_{jt}\},~\forall j,~\forall t \\
\hat w_{jt}=W_{jt}^e+w_{jt}^u-w_{jt}^l,\forall j,~\forall t \\
0 \le w_{jt}^l, w_{jt}^u \le W_{jt}^h,\forall j,\forall t \\
w_{jt}^lw_{jt}^u=0,~\forall j,~\forall t\\
\sum_{j=1}^{N_w} \frac{w_{jt}^l+w_{jt}^u}{W_{jt}^h} \le \Gamma_S, ~\forall t\\
\sum_{t=1}^{T} \frac{w_{jt}^l+w_{jt}^u}{W_{jt}^h} \le \Gamma_T, \forall j\\
\end{array}} \right.} \right\}
\end{align}
where $W_{jt}^e = (W_{jt}^l+W_{jt}^u)/2$, $W_{jt}^h=(W_{jt}^u-W_{jt}^l)/2,\forall j,\forall t$.

The RG output after preparatory curtailment is the minimum of the actual RG output and the upper limit, as in the first constraint. The remaining constraints limit the deviation of the actual RG output that fluctuates within a range $[W_{jt}^l, W_{jt}^u]$, whose mean is the forecast value $W_{jt}^e$. When the RG output is higher than its forecast, $w_{jt}^u>0, w_{jt}^l=0$; otherwise, $w_{jt}^u=0, w_{jt}^l>0$. To
avoid over-conservativeness of the model, uncertainty budgets $\Gamma_S$ and $\Gamma_T$ are added to restrain the spatial and temporal deviations from the prediction.

The above uncertainty set involves nonlinear constraints $w_{jt}=\mbox{min}\{\hat w_{jt},\xi_{jt}\},\forall j,\forall t$, each of which can be linearized by introducing a binary variable as below.
\bsq
\label{eq:uncertainty}
\begin{align}
    w_{jt} \le \hat w_{jt}, w_{jt} \le \xi_{jt} \label{eq:uncertainty-1}\\
    \hat w_{jt} -(W_{jt}^u-W_{jt}^l)(1-z_{jt}) \le w_{jt} \label{eq:uncertainty-2}\\
    \xi_{jt} - (\xi_{jt}-W_{jt}^l) z_{jt} \le w_{jt} \label{eq:uncertainty-3}\\
    z_{jt} \in \{0,1\}
\end{align}
\esq
when $z_{jt}=0$, constraints \eqref{eq:uncertainty-1} and \eqref{eq:uncertainty-3} lead to $w_{jt}=\xi_{jt}\le \hat w_{jt}$, while \eqref{eq:uncertainty-2} becomes $ \hat w_{jt}-w_{jt} \le W_{jt}^u-W_{jt}^l$ which is always satisfied. In this case, there is surplus RG output and curtailment is executed. When $z_{jt}=1$, constraints \eqref{eq:uncertainty-1} and \eqref{eq:uncertainty-2} lead to $w_{jt}=\hat w_{jt} \le \xi_{jt}$, while \eqref{eq:uncertainty-3} becomes $w_{jt} \ge W_{jt}^l$ which is always satisfied. Thus, no RG output is curtailed. The other nonlinear constraint $w_{jt}^lw_{jt}^u=0$ can be linearized as
\begin{align}
    0 \le w_{jt}^l \le M \hat z_{jt},~ 0 \le w_{jt}^u \le M (1-\hat z_{jt}), \hat z_{jt} \in \{0,1\}
\end{align}
where $M$ is a large constant.

After the above transformation, the uncertainty set \eqref{eq:uncertainty} consists of mixed-integer linear constraints. Moreover, it relies on the value of the first-stage curtailment strategy $\xi$, and thus, is \emph{decision-dependent}. Traditional algorithms for solving robust optimizations with decision-independent uncertainty set are the Bender’s decomposition \cite{bertsimas2012adaptive} and the C\&CG algorithm \cite{zeng2013solving}. However, they cannot be directly applied to solve the robust optimization with DDU due to the following two difficulties:

1) \emph{Failure in robust feasibility and optimality}. When the first-stage decision changes, the previously selected scenarios may be outside of the new uncertainty set. For example, suppose $w_1$ is the worst-case scenario corresponding to a first-stage decision $\xi_1$, i.e., $w_1 \in \mathcal{W}(\xi_1)$. Then when the decision changes to $\xi_2$, we could have $w_1 \notin \mathcal{W}(\xi_2)$. Consequently, adding the scenario $w_1$ to the master problem may no longer result in a relaxation of the original robust optimization. This could lead to over-conservative or even infeasible solutions.

2) \emph{Failure in finite convergence}. The proof of convergence of the C\&CG algorithm is based on the fact that: the worst scenario always rests at a vertex of the uncertainty set, the same vertex will not be picked up twice, and the number of vertices is finite. However, when the uncertainty set becomes decision-dependent, the previously selected scenarios may no longer be at the vertex of the new set, and so there is no convergence guarantee.

In this paper, we develop an adaptive C\&CG algorithm to address the above issues, which will be explained in detail in Section \ref{sec-3}.

\subsection{Robust generation dispatch model}
 The robust generation dispatch model considering preparatory curtailment of RGs is given below: \footnote{$[A]$ denotes all positive integers that do not exceed $A$.}
\bsq
\label{eq:ED}
\begin{align}
    \mathop{\min}_{p, r^{\pm}, \xi}~ & \sum_{t=1}^T \sum_{i=1}^{N_g} \left(\alpha_ip_{it}+\beta_i^{+}r_{it}^{+}+\beta_i^{-}r_{it}^{-}\right)+\sum_{t=1}^T \sum_{j=1}^{N_w} \gamma_j(W_{jt}^u-\xi_{jt})^2   \nonumber \\
    ~ & + \mathop{\max}_{w \in \mathcal{W}(\xi)} Q(p,r^{\pm},w) \label{eq:ED-1}\\
    \mbox{s.t.}~ & W_{jt}^l \le \xi_{jt} \le W_{jt}^u, ~\forall j \in [N_w],~\forall t \in [T] \label{eq:ED-2}\\
    ~ & u_{it}P_i^l \le p_{it} \le u_{it}P_i^u,~\forall i \in [N_g],~\forall t \in [T] \label{eq:ED-3}\\
    ~ & \left(p_{it}+r_{it}^{+}\right)-\left(p_{i(t-1)}-r_{i(t-1)}^{-}\right) \le u_{i(t-1)}\mathcal{R}_i^u  \nonumber\\
    ~ & +\left(1-u_{i(t-1)}\right)P_i^u,\forall i \in [N_g], \forall t \in [T]/\{1\} \label{eq:ED-3-2.1}\\
    ~ & -\left(p_{it}+r_{it}^{+}\right)+\left(p_{i(t-1)}-r_{i(t-1)}^{-}\right) \le u_{it}\mathcal{R}_i^d  \nonumber\\
    ~ & +\left(1-u_{it}\right)P_i^u,\forall i \in [N_g], \forall t \in [T]/\{1\} \label{eq:ED-3-2.2}\\
        ~ & \left(p_{it}-r_{it}^{-}\right)-\left(p_{i(t-1)}+r_{i(t-1)}^{+}\right) \le u_{i(t-1)}\mathcal{R}_i^u \nonumber\\
    ~ & +\left(1-u_{i(t-1)}\right)P_i^u,\forall i \in [N_g], \forall t \in [T]/\{1\}\label{eq:ED-3-2.3}\\
    ~ & -\left(p_{it}-r_{it}^{-}\right)+\left(p_{i(t-1)}+r_{i(t-1)}^{+}\right) \le u_{it}\mathcal{R}_i^d  \nonumber\\
    ~ & +\left(1-u_{it}\right)P_i^u,\forall i \in [N_g], \forall t \in [T]/\{1\}\label{eq:ED-3-2.4}\\
    ~ & u_{it}P_i^l \!\!\le\! p_{it}-r_{it}^{-},\! p_{it}+r_{it}^{+} \!\!\le\! u_{it}P_i^u,\!\forall i \in [N_g],\forall t \in [T] \label{eq:ED-4}\\
    ~ & 0 \!\le\! r_{it}^{-}\! \le\! u_{it}R_{i}^{-}\!,\!0\! \le\! r_{it}^{+} \!\le\! u_{it}R_i^{+},\!\forall i \in [N_g],\forall t \in [T] \label{eq:ED-5}\\
    ~ & \sum \nolimits_{j=1}^{N_w} W_{jt}^e + \sum \nolimits_{i=1}^{N_g} p_{it} =\sum \nolimits_{l=1}^{N_l} \mathcal{D}_{lt},~\forall t \in [T] \label{eq:ED-6}\\
    ~ & -F_k \le \sum \nolimits_{i=1}^{N_g} \pi_{ik} p_{it} +\sum \nolimits_{j=1}^{N_w}\pi_{jk} W_{jt}^e \nonumber\\ ~ & -\sum \nolimits_{l=1}^{N_l}\pi_{lk} \mathcal{D}_{lt} \le F_k,~\forall k \in [N_k],~\forall t \in [T]\label{eq:ED-7}
\end{align}
\esq
The first-stage decisions, including the contemporary output $\{p_{it},\forall i,\forall t\}$ and reserve capacity $\{r_{it}^{\pm},\forall i,\forall t\}$ of thermal units, and the preparatory curtailment schedule $\{\xi_{jt},\forall j,\forall t\}$, are made prior to the realization of uncertainty. The $\{u_{it},\forall i,\forall t\}$  are binary parameters and $u_{it}=1/0$ indicates the generator is on/off. Here, the unit commitment $\{u_{it},\forall i,\forall t\}$ are given and fixed parameters. The objective function \eqref{eq:ED-1} minimizes the sum of generation cost, reserve cost, preparatory curtailment cost, and the re-dispatch generation and curtailment costs under the worst-case scenario. The $\gamma_j$ is a coefficient that fully incorporates the opportunity cost of renewable power curtailment, including but not limited to the resulting increased carbon emission, the impact on power system stability, and the possibility of exceeding the rate of renewable power spilling limits. A quadratic function that is widely used in literature, such as \cite{li2015adjustable}, is adopted. Denote the objective of the first-stage problem as $f_s$, then given the $\gamma_j>0,\forall j$, we have
\begin{align}
   \frac{\partial f_s}{\partial \xi_{jt}}=-\gamma_j(W_{jt}^u-\xi_{jt})\le 0 < \frac{\partial f_s}{\partial p_{it}}=\alpha_i \nonumber
\end{align}
According to equal incremental principle, it assures that renewable generators have higher priority compared with the thermal generators.

The curtailment strategy is between the lower and upper bounds of the actual RG output as in constraint \eqref{eq:ED-2}. Constraints \eqref{eq:ED-3} and \eqref{eq:ED-5} represent the limits of generation capacity and reserve capacity, respectively. Constraints \eqref{eq:ED-3-2.1}-\eqref{eq:ED-3-2.4} ensure that the ramping constraints still hold when offering the reserve capacities in the second-stage, so we do not need to consider the ramping rates in the second stage. The generation adequacy is ensured in \eqref{eq:ED-4}. Power balancing condition and power flow limits are given in \eqref{eq:ED-6} and \eqref{eq:ED-7}, respectively. The term
$Q(p,r^{\pm},w)$ in \eqref{eq:ED-1} is the re-dispatch cost when the RG output is $w$. Therefore, maximizing $Q(p,r^{\pm},w)$ over all possible $w \in \mathcal{W}(\xi)$ gives the cost under the worst-case scenario. Here, $\mathcal{W}(\xi)$ is the decision-dependent uncertainty set \eqref{eq:DDU-set}. To be specific, $Q(p,r^{\pm},w)$ is the optimal objective value of the second-stage re-dispatch problem:
\bsq
\label{eq:SD}
\begin{align}
    \mathop{\min}_{p^{\pm},\hat \xi}~ & \sum_{t=1}^T \sum_{i=1}^{N_g} \left(d_i^{+} p_{it}^{+} +d_i^{-}p_{it}^{-}\right) + \sum_{t=1}^T \sum_{j=1}^{N_w} \hat \gamma_j \hat \xi_{jt}^2 \label{eq:SD-1}\\
    \mbox{s.t.}~ & 0 \le p_{it}^{+} \le r_{it}^{+},~0 \le p_{it}^{-} \le r_{it}^{-},~\forall i \in [N_g],~\forall t \in [T]\label{eq:SD-2}\\
    ~ & \hat \xi_{jt} \ge 0,~\forall i \in [N_w],~\forall t \in [T] \label{eq:SD-2-2}\\
    ~ & \sum_{i=1}^{N_g} \left(p_{it}+p_{it}^{+}-p_{it}^{-}\right) + \sum_{j=1}^{N_w} (w_{jt}-\hat \xi_{jt}) = \sum_{l=1}^{N_l} \mathcal{D}_{lt},\forall t \in [T] \label{eq:SD-3}\\
    ~ & -F_k \le \sum_{i=1}^{N_g} \pi_{ik} \left(p_{it}+p_{it}^{+}-p_{it}^{-}\right)  +\sum_{j=1}^{N_w} \pi_{jk} (w_{jt}-\hat \xi_{jt}) \nonumber\\
    ~ & - \sum_{l=1}^{N_l} \pi_{lk}\mathcal{D}_{lt}  \le F_k,~\forall k \in [N_k],~\forall t \in [T] \label{eq:SD-4}
\end{align}
\esq
The second-stage decisions, including the incremental outputs of thermal units $\{p_{it}^{\pm},\forall i,\forall t\}$ and real-time curtailment $\{\hat \xi_{jt},\forall j,\forall t\}$, are made after knowing the exact outputs of RGs. The objective function \eqref{eq:SD-1} is to minimize the the re-dispatch
cost, i.e., the sum of up- and down-regulation costs and real-time curtailment costs. Constraint \eqref{eq:SD-2} requires the incremental output of thermal unit be within its reserve capacity. Power balancing condition and power flow limits are \eqref{eq:SD-3} and \eqref{eq:SD-4}. There is a quadratic term $\hat \gamma_j \hat \xi_{jt}^2$ in the objective function, which can be linearized via a convex combination approach \cite{wu2011tighter}. To be specific, suppose we have evaluations at $N$ discrete points $\xi_{jt}^1,...,\xi_{jt}^{N}$, and $\hat g_{jt}^1= \hat \gamma_j (\hat \xi_{jt}^1)^2$,···,$\hat g_{jt}^{N}= \hat \gamma_j (\hat \xi_{jt}^{N})^2$. By introducing variables $\sigma_{jt}^1,...,\sigma_{jt}^{N} \ge 0$ 
and $\sum_{n=1}^{N} \sigma_{jt}^n=1$, $\hat \xi_{jt}$ and $\hat \gamma_j \hat \xi_{jt}^2$ can
be replaced with linear functions $\sum_{n=1}^{N} \sigma_{jt}^n \xi_{jt}^n$ and $\sum_{n=1}^{N} \sigma_{jt}^n g_{jt}^n$
in $\sigma$, respectively.


\subsection{Compact form}
The aforementioned robust generation dispatch model with decision-dependent uncertainty can be written in the following compact form:
\begin{align}
\label{eq:robust-whole}
    \min_{x \in \mathbb{X}, \xi \in \Xi} c^\top x+g(\xi) + \max_{w \in \mathcal{W}(\xi)} \min_{y \in \mathbb{F}(x,w)} d^\top y
\end{align}
and 
\begin{align}
\mathbb{X}=~ & \{x \in \mathbb{R}^{3N_g \times T}~:~ Ax\ge b\} \nonumber\\
\Xi=~ & \{\xi \in \mathbb{R}^{N_w\times T}~:~ B\xi \le e\} \nonumber\\
\mathcal{W}(\xi)=~ & \{w \in \mathbb{R}^{N_w\times T}~:~ Hw \le C \hat w + F\xi + G z+a, \nonumber\\
~ &  \forall z \in \mathbb{Z}^{2N_w\times T}, ~\hat w \in \mathbb{R}^{3N_w\times T}\} \nonumber\\
\mathbb{F}(x,w)=~ & \{y \in \mathbb{R}^{(2N_g+N\times N_w)\times T}~:~Ey \le f-Rw-Dx\}
\end{align}
where $x:=\{p_{it},r_{it}^{+}, r_{it}^{-},\forall i,\forall t\}$ and $c$ is a collection of $\alpha_i$,$\beta_i^{+}$,$\beta_i^{-},\forall i$. $\xi$ is a collection of $\xi_{jt},\forall j,\forall t$; $g(\xi)=(W^u-\xi)^\top \textbf{diag}(\gamma_1,...,\gamma_{N_w})(W^u-\xi)$. $w:=\{w_{jt},\forall j,\forall t\}$, $\hat w:=\{\hat w_{jt},w_{jt}^l,w_{jt}^u,\forall j,\forall t\}$, and $z$ is a collection of the binary variables $z_{jt}$ and $\hat z_{jt}$ for all $j,t$. $y:=\{p_{it}^{+},p_{it}^{-},\forall i,\forall t; \sigma_{jt}^n,\forall j,\forall n,\forall t\}$; $d$ summarizes $d_i^{+},d_i^{-},\forall i$ and $g_{jt}^n,\forall j,\forall n,\forall t$. $A,B,H,C,F,G,E,R,D$ and $d, b,e,a,f$ are the coefficient matrices and vectors.

\textbf{Remark:} It is worth noting that the proposed model \eqref{eq:robust-whole} subsumes the two-stage robust optimization with decision-independent uncertainty (DIU) as a special case. When $F$ in $\mathcal{W}(\xi)$ is an all-zero matrix, the DDU set degenerates to a DIU set. The algorithm proposed in this paper can still be applied.

\section{Solution Algorithm}
\label{sec-3}
Different from the traditional robust optimization models \cite{wei2014robust}, the uncertainty set in \eqref{eq:robust-whole} is \emph{decision-dependent}. The widely-used algorithms for solving the traditional robust optimizations may cause oscillation because the update of the first-stage decision keeps changing the shape as well as the extreme points of the uncertainty set. In this section, an adaptive C\&CG algorithm is proposed to generate the optimal solution in a finite number of iterations.

In this paper, we do not require relatively complete recourse, which assumes that the second-stage problem \eqref{eq:SD} is feasible for any fixed first-stage action $x$ and scenario $w$. Instead, we assume that

\noindent \textbf{A1}: The robust optimization problem \eqref{eq:robust-whole} has at least one \emph{robust feasible} solution, i.e., there exists a $x \in \mathbb{X}, \xi \in \Xi$ such that for any $w \in \mathcal{W}(\xi)$, the second-stage problem \eqref{eq:SD} is always feasible.

This condition better matches engineering practices.

\subsection{Equivalent transformation}
Given the first-stage decisions $x \in \mathbb{X}$ and $\xi \in \Xi$, the inner ``max-min'' problem is a bilevel optimization:
\begin{align}
\label{eq:max-min}
    S(x,\xi):=\max_{w \in \mathcal{W}(\xi)} \min_{y \in \mathbb{F}(x,w)} d^\top y 
\end{align}
where $\mathbb{F}(x,w)$ is a polyhedral set while $\mathcal{W}(\xi)$ is bounded and consists of mixed-integer linear constraints depending on the value of $\xi$ in the first-stage.

Note that the first-stage decisions $x \in \mathbb{X}$ and $\xi \in \Xi$ may not always be robust feasible. Therefore, we construct the following problem for feasibility-check first.
\begin{align}\label{eq:feasibility-1}
    \min_{s} ~ & 1^{\top}s \nonumber\\
    \mbox{s.t.}~ & Ey \le f-Rw-Dx+s,~ s \ge 0
\end{align}
where $s$ is the slack variable. Obviously, the problem \eqref{eq:feasibility-1} is always feasible. If the original second-stage problem \eqref{eq:SD} is feasible, we have $s$ is an all-zero vector with a zero optimal objective value; otherwise, the optimal objective value is larger than zero. Then the inner “max-min” problem for feasibility-check is equivalent to
\begin{align}\label{eq:feasibility-2}
    S_f(x,\xi):=\!\!\max_{w \in \mathcal{W}(\xi), y,s, \atop \mu, \mu_s, z_{\mu}, z_{\mu_s}} & 1^{\top}s \nonumber\\
    \mbox{s.t.}~ & E^{\top} \mu=0,~ 1-\mu-\mu_s=0 \nonumber\\
    ~ & 0 \le \mu \le \hat Mz_{\mu},~ 0 \le \mu_s \le \hat M z_{\mu_s} \nonumber\\
    ~ & 0 \le (-Ey+f-Rw-Dx+s) \le \hat M (1-z_{\mu}) \nonumber\\
    ~ & 0 \le s \le \hat M (1-z_{\mu_s}) \nonumber\\
    ~ & z_{\mu}, z_{\mu_s} \in \{0,1\}^{n_{\mu}}
\end{align}
where $n_{\mu}$ is the number of constraints in $\mathbb{F}(x,w)$, $\mu, \mu_s$ are dual variables, and $\hat M$ is a large enough constant. Problem \eqref{eq:feasibility-2} is obtained by replacing the problem \eqref{eq:feasibility-1} with its KKT condition and linearize the complementary and slackness condition in the form of $0 \le x \perp y \ge 0$ using the Big-M method \cite{fortuny1981representation}. Then $S_f(x,\xi)$=0 if and only if $(x,\xi)$ is robust feasible. Denote the above optimization problem \eqref{eq:feasibility-2} as the feasibility-check problem \textbf{FCP}.

If $(x,\xi)$ is robust feasible, we can then solve the original ``max-min'' problem \eqref{eq:max-min}, which is equivalent to
\bsq
\label{eq:subproblem}
\begin{align}
    \max_{w \in \mathcal{W}(\xi), y, \mu, z_{\mu}} ~ &  d^\top y \label{eq:subproblem-1}\\
    \mbox{s.t.} ~ & d +  E^\top \mu =0 \label{eq:subproblem-2}\\
    ~ & 0 \le \mu \le \hat M z_{\mu} \label{eq:subproblem-3}\\
    ~ & 0 \le (-Ey+f-Rw-Dx) \le \hat M (1-z_{\mu}) \label{eq:subproblem-4}\\
    ~ & z_{\mu} \in \{0,1\}^{n_\mu} \label{eq:subproblem-5}
\end{align}
\esq
We denote the above optimization \eqref{eq:subproblem} as the subproblem \textbf{SP}.

\textbf{Remark on the selection of the Big-M}: The big-M method is applied for linearizing the constraint $w_{jt}^l w_{jt}^u=0$ and also the complementary slackness condition in KKT. For the $w_{jt}^l w_{jt}^u=0$, since both $w_{jt}^l$ and $w_{jt}^u$ are in $[W_{jt}^l,W_{jt}^u]$, we can just let $M=W_{jt}^u$. The complementary slackness condition in KKT is linearized using \eqref{eq:subproblem-3}-\eqref{eq:subproblem-5}. The constraint $-Ey+f-Rw-Dx \ge 0$ refers to the physical constraints in problem \eqref{eq:SD}, whose upper bounds can be derived according to their physical meanings. For example, for the constraint $p_{it}^{+}-r_{it}^{+} \le 0$, we have $0 \le r_{it}^{+}-p_{it}^{+} \le r_{it}^{+} \le R_{it}^{+}$. The variable $\mu$ is the dual variable of the physical constraints, it reflects the marginal effect on the primal objective \eqref{eq:SD-1} per unit change in the primal constraint limit. Therefore, the upper bound of $\mu$ can be obtained according to some sensitivity analysis on the power system beforehand. The big-M can be chosen as the maximum of the upper bounds for $-Ey+f-Rw-Dx$ and for $\mu$. It’s worth noting that a larger $M$ is not always the better, since it may cause numerical instabilities. In the case study, we let $M=W_{jt}^u$ with $\hat M=10^5$ for the IEEE 39-bus system and $\hat M=10^6$ for the IEEE 118-bus system. All simulations can output the correct solutions efficiently.

\begin{proposition}
\label{prop-1}
The optimal solution $w^*$ of the problem \textbf{FCP} or \textbf{SP} will be reached at one of the extreme points of $\mathcal{W}(\xi)$.
\end{proposition}
\begin{proof}
Given $x$ and $\xi$, the uncertainty set $\mathcal{W}(\xi)$ can be seen as the union of $2^{2N_w \times T}$ polytopes, each of which corresponds to a particular assignment to the set of binary variables $z \in \mathbb{Z}^{2N_w \times T}$, i.e., $\mathcal{W}(\xi)=\mathcal{W}^1(\xi) \cup \mathcal{W}^2(\xi) \cup ... \cup \mathcal{W}^{2^{2N_w\times T}} (\xi)$. For the subproblem \textbf{SP}, $Q(x,w):=\min_{y \in \mathbb{F}(x,w)} d^\top y$ is a quasiconvex function in $w$, so the optimal solution of $\max_{w \in \mathcal{W}^p(\xi)} Q(x,w)$ will be achieved at the extreme point of $W^p(\xi)$ \cite{rockafellar2015convex}. The $\mbox{argmin}_{w \in \mathcal{W}^p(\xi)} Q(x,w),\forall p\in [2^{2N_w\times T}]$ with the highest objective value is the optimal solution of \eqref{eq:max-min}, which is also an extreme point of $\mathcal{W}(\xi)$. Similarly, we can prove the conclusion for the feasibility-check problem \textbf{FCP}.
\end{proof}

Suppose the worst-case scenario identified by the \textbf{FCP} or the \textbf{SP} is $w^*$ and the corresponding value of the binary variable is $z^*$. At the optimal point, the set of active and inactive constraints in $\mathcal{W}(\xi)$ are
\begin{align}
    H_{eq} w = ~ &  C_{eq} \hat w + F_{eq}\xi + G_{eq} z+ a_{eq} \label{eq:active}\\
    H_{in} w < ~ &  C_{in} \hat w + F_{in}\xi + G_{in} z+ a_{in} \label{eq:inactive}
\end{align}

If we fix the value of $z$ to $z^*$, then in non-degenerate cases, the matrix $[H_{eq},-C_{eq}]$ is a ($4N_w \times T$)-dimensional full rank matrix as the active constraints form a set of basis. For any given $\xi$, the value of $(w, \hat w)$ can be uniquely determined by
\begin{align}
\left[\begin{array}{c}
w\\
\hat w
\end{array} \right]=[H_{eq},-C_{eq}]^{-1} (F_{eq}\xi + G_{eq}z^* +a_{eq})
\end{align}

However, in degenerate cases, it could be hard to recover $w$ and $\hat w$ using the active constraints. Methods to deal with this challenge are introduced in the following. Denote the set of active constraints we finally obtain in iteration $n$ as $\mathcal{AC}_n$.

\subsection{Model degeneracy}
\label{sec-3-2}
\begin{figure}[t]
\centering
\includegraphics[width=1.0\columnwidth]{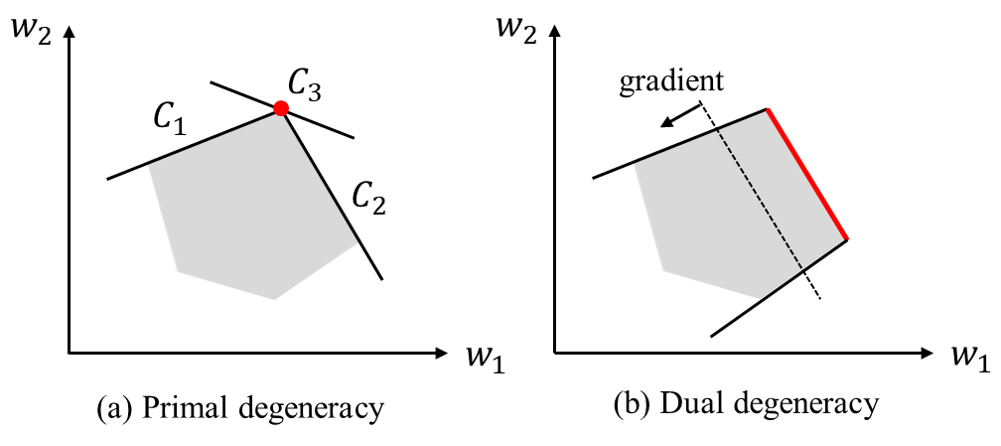}
\caption{Illustration of degenerate cases.}
\label{fig:degeneracy}
\end{figure}
In general, two types of degeneracy may occur, i.e., primal degeneracy and dual degeneracy, which are illustrated in Fig. \ref{fig:degeneracy}. Primal degeneracy is caused by weak redundant constraints. For example, in Fig. \ref{fig:degeneracy}(a), constraints C1-C3 are all active at the optimal point so equation \eqref{eq:active} is over-determined. In fact, constraint C3 is redundant: although C3 intersects the feasible region, removing C3 will not change the feasible region. The impact of primal degeneracy is manageable, because redundant constraints in $\mathcal{W}(\xi)$ can be filtered before identifying the active constraints, and then $[H_{eq},-C_{eq}]$ is invertible. The method for removing redundant constraints is given in Appendix A. In particular, we need to be careful about three cases: 
\begin{enumerate}
    \item If $w_{jt}^{*}=W_{jt}^u$ : In this case, we must have $w_{jt}^{*}=\xi_{jt}$; otherwise, we have $w_{jt}^{*}=\hat w_{jt}^{*} \le \xi_{jt} \le W_{jt}^u$. Equality holds for all inequalities in the middle, so we still have $w_{jt}^{*}=\xi_{jt}=W_{jt}^u$. Therefore, we have three active constraints $w_{jt}=\hat w_{jt}$ and $\hat w_{jt}=W_{jt}^e+w_{jt}^u-w_{jt}^l$, and $w_{jt}^l=0$. Moreover, the active constraints $w_{jt}^u=W_{jt}^h$, and $w_{jt}=\xi_{jt}$ give the same limits and one of them is redundant. Our goal is to add the active constraint \eqref{eq:active} to the master problem to ensure that when the first-stage decision changes, the scenario generated remains at the extreme point of the new set. To this end, here we remove the $w_{jt}^u=W_{jt}^h$. 
    \item If $\xi_{jt}=W_{jt}^l$ and $w_{jt}^{l*}=W_{jt}^h$: Since $\hat w_{jt}=W_{jt}^e+w_{jt}^u-w_{jt}^l$ and $w_{jt}^uw_{jt}^l=0$, we have $\hat w_{jt}^*=W_{jt}^l=\xi_{jt}$. Therefore, we have three active constraints $w_{jt}=\hat w_{jt}$, $\hat w_{jt}=W_{jt}^e+w_{jt}^u-w_{jt}^l$, and $w_{jt}^u=0$. Moreover, the active constraints $w_{jt}^l=W_{jt}^h$ and $w_{jt}=\xi_{jt}$ give the same limits. For the reasons in 1), we remove the $w_{jt}^l=W_{jt}^h$.
    \item If $\xi_{jt}=w_{jt}^e$ with $w_{jt}^{l*}=w_{jt}^{u*}=0$ : Since $\hat w_{jt}=W_{jt}^e+w_{jt}^u-w_{jt}^l$, we have $\hat w_{jt}^{*}=W_{jt}^e$. The two active constraints $w_{jt}^l=0$ and $w_{jt}^u=0$ gave the same limit as $w_{jt}=\xi_{jt}$. For the reasons in 1), we keep the $w_{jt}=\xi_{jt}$. For the constraint $w_{jt}^l=w_{jt}^u=0$, due to the complementary condition $w_{jt}^lw_{jt}^u=0$, we keep one of them and remove the other.
\end{enumerate}


Dual degeneracy occurs if multiple solutions attain the optimal value, then the optimal
solution is not unique. As in Fig. \ref{fig:degeneracy}(b), a facet is perpendicular to the gradient of the objective function, so any point on the facet is an optimal solution. When there is dual degeneracy, the matrix $[H_{eq},-C_{eq}]$ is not full-rank, and \eqref{eq:active} is underdetermined. We can include some of the constraints in \eqref{eq:inactive} to make $[H_{eq},-C_{eq}]$ a full rank matrix, i.e., letting $H_{in,j}w=C_{in,j}\hat w + F_{in,j} \xi +G_{in,j} z^* + a_{in,j}$ for some $j$. If the unique $w^*$ and $\hat w^*$ we get satisfy constraint \eqref{eq:inactive}, then use the new set of constraints as the active constraint.

\subsection{Master Problem}
The master problem (\textbf{MP}) in iteration $N$ is formulated as
\bsq
\label{eq:master}
\begin{align}
    \min_{x, \xi} ~ &  c^\top x +g(\xi) +\eta \label{eq:master-1}\\
    \mbox{s.t.} ~ &  x \in \mathbb{X}, \xi \in \Xi \label{eq:master-2}\\
    ~ & \eta \ge d^\top y^n, \forall n \in [N] \label{eq:master-3}\\
    ~ & H_{eq}^n \tilde w^n = C_{eq}^n \hat w^n + F_{eq}^n \xi + G_{eq}^nz^{n*} +a_{eq}^n ,\forall n \in [N] \label{eq:master-4}\\
    ~ & w^n = \min\{\tilde w^n, \xi\}, w^n =\max\{W^l, \tilde w^n\},\forall n \in [N] \label{eq:master-5}\\
    ~ & y^n \in \mathbb{F}(x, w^n), \forall n \in [N] \label{eq:master-6}
\end{align}
\esq
The decision variables are $\{x,\xi,\eta,y^n,w^n,\hat w^n,\tilde w^n, \forall n\in [N]\}$. $n$ is the index of scenarios. $w^n$ is the RG output scenario in the $n$-th iteration while $\tilde w^n$ and $\hat w^n$ are intermediate variables. The constraints of the first-stage problem are given in \eqref{eq:master-2}. Since $\eta$ is minimized in the objective \eqref{eq:master-1}, constraint \eqref{eq:master-3} is equivalent to $\eta=\max_n \min_{y^n} d^\top y^n$ which is the re-dispatch cost under the worst-case scenario. The constraints in the second-stage problem are given in \eqref{eq:master-6}. Our key innovation lies in the scenario representation \eqref{eq:master-4}-\eqref{eq:master-5}. Instead of adding the worst-case scenario in the previous iteration itself directly to the master problem as the traditional robust optimization methods \cite{zeng2013solving}, here we include the set of active constraints \eqref{eq:master-4}. When the first-stage decision $\xi$ changes, the $\tilde w^n$ will change accordingly and can be uniquely determined by \eqref{eq:master-4}. Together with \eqref{eq:master-5}, we can ensure that the scenario $w^n$ is still an extreme point of $\mathcal{W}(\xi)$. We will give an example to explain this in detail in the Remark later.

The constraint $w^n=\min\{\tilde w^n, \xi\}$ can be linearized as
\bsq
\label{eq:linearize2}
\begin{align}
    W^l \le \tilde w^n \le W^u, ~z_w^n \in \{0,1\}^{N_w \times T} \label{eq:linearize2-1}\\
    w^n \le \tilde w^n,~ w^n \le \xi \label{eq:linearize2-2}\\
    \tilde w^n -(W^u-W^l)z_w^n \le w^n \label{eq:linearize2-3}\\
    \xi -(W^u-W^l)(1-z_w^n) \le w^n  \label{eq:linearize2-4}
\end{align}
\esq
When $\tilde w^n \le \xi$, let $z_w^n=0$, constraint \eqref{eq:linearize2-3} turns into $\tilde w^n \le w^n$. Together with \eqref{eq:linearize2-2}, we have $w^n=\tilde w^n=\min\{\tilde w^n, \xi\}$. Moreover, constraint \eqref{eq:linearize2-4} becomes $\xi-(W^u-W^l) \le \tilde w^n$, which is satisfied since both $\xi$ and $\tilde w^n$ are in $[W^l, W^u]$.
When $\tilde w^n \ge \xi$, let $z_w^n=1$, constraint \eqref{eq:linearize2-4} turns into $\xi \le w^n$. Together with \eqref{eq:linearize2-2}, we have $w^n=\xi=\min\{\tilde w^n, \xi\}$. Moreover, constraint \eqref{eq:linearize2-3} becomes $\tilde w^n-(W^u-W^l) \le \xi$, which is satisfied since both $\xi$ and $\tilde w^n$ are in $[W^l, W^u]$.

This is different from the linearization technique in \eqref{eq:uncertainty} since now $\xi$ is a variable instead of a constant. Constraint $w^n=\max\{W^l, \tilde w^n\}$ can be linearized in a similar way since it is equivalent to $-w^n=\min\{-W^l, -\tilde w^n\}$. The master problem \textbf{MP} is thus turned into a mixed-integer quadratic program (MIQP) that can be solved efficiently by commercial software.

\begin{figure*}[t]
\centering
\includegraphics[width=1.8\columnwidth]{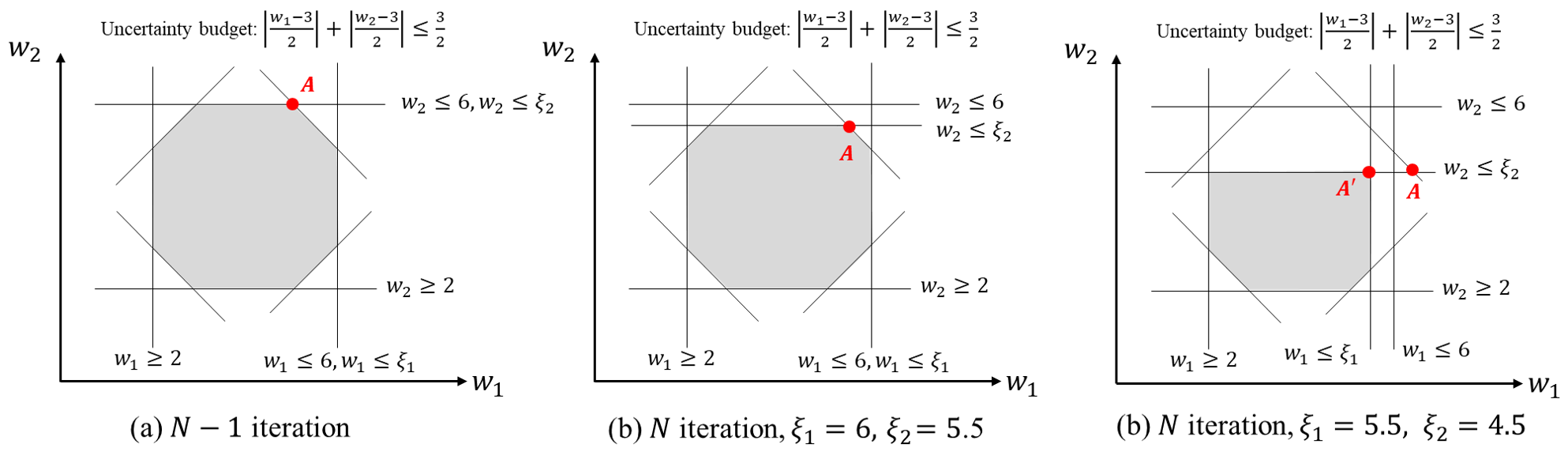}
\caption{Illustration of the adaptive RG power output scenario.}
\label{fig:uncertainty}
\end{figure*}
\textbf{Remark:} Since $[H_{eq}^n,-C_{eq}^n]$ is a full-rank matrix, $\tilde w^n$ and $\hat w^n$ can be uniquely represented as a function of $\xi$. If $\tilde w^n$ is a feasible point of the new $\mathcal{W}(\xi)$, then we have $\tilde w^n \le \xi$ and $\tilde w^n \ge W^l$ due to the constraints in $\mathcal{W}(\xi)$. Therefore, $w^n=\tilde w^n$ is an extreme point of the new $\mathcal{W}(\xi)$. If $\tilde w^n$ is infeasible for the new $\mathcal{W}(\xi)$, then constraint \eqref{eq:master-5} projects it onto an extreme point (vertex) of the new $\mathcal{W}(\xi)$. We use a simple example with two RGs for illustration:


Suppose $T=1$ and the uncertainty set is \footnote{The $t$ index is omitted for simplicity.}
\begin{align}
    \mathcal{W}(\xi)=\left\{ {(w_1,w_2)^\top ~\left| {\begin{array}{*{20}{c}}
     w_j = \min\{\hat w_j, \xi_j\},\forall j=1,2 \\
     \hat w_{j} = 4 + w_{j}^u-w_j^l,\forall j=1,2 \\
   0 \le w_j^u, w_j^l \le 2,\forall j=1,2\\
   \frac{w_1^l+w_1^u}{2}+\frac{w_2^l+w_2^u}{2} \le \frac{3}{2} \\
\end{array}} \right.} \right\}
\end{align}

Suppose in the $N-1$ iteration, $\xi_1=\xi_2=6$, so the feasible region is the grey region in Fig. \ref{fig:uncertainty}(a). Suppose the optimal solution of the \textbf{SP} is point A. Primal degeneracy happens in this case and the constraints $w_2 \le \xi_2$ and $w_2 \le 6$ give the same limit. According to Section \ref{sec-3-2}, we keep the former constraint and drop the latter one. Therefore, the scenario $A$ can be represented as a function of $\xi$, i.e. $(\tilde w_1,\tilde w_2)=(11-\xi_2, \xi_2)$. Then in the $N$ iteration, if $\xi=(6,5.5)$ as in Fig. \ref{fig:uncertainty}(b), the active constraints remain unchanged. Substituting the value of $\xi$ into the expression of point A, we can get a new point A $(\tilde w_1,\tilde w_2)=(5.5,5.5)$, which is a vertex of the new $\mathcal{W}(\xi)$. If $\xi=(5.5,4.5)$ as in Fig. \ref{fig:uncertainty}(c), now the new point A is $(\tilde w_1,\tilde w_2)=(6.5, 4.5)$ which is outside of $\mathcal{W}(\xi)$. But with constraint \eqref{eq:master-5}, the point A ($\tilde w$) is projected to point A' ($w$) which is again a vertex of the new set $\mathcal{W}(\xi)$. 

\subsection{Adaptive C\&CG algorithm}
\label{sec-3-4}
The overall iterative algorithm is given in Algorithm 1.
\begin{algorithm}
\normalsize
\label{Ag:Path-Generation}
\caption{{\bf : Adaptive C\&CG Algorithm}}
\begin{algorithmic}[1]
\STATE \textbf{Initiation}:  Error tolerance $\varepsilon>0$; $N=1$; Let $C_{eq}^N$, $F_{eq}^N$, $G_{eq}^N$ be zero matrices; Let $H_{eq}^N$ be an identity matrix and $a_{eq}^N=W^e$; Assign a large value to $UB_0$.
 
\STATE \textbf{Solve the Master Problem}

Solve the $\textbf{MP}$ \eqref{eq:master}. Let $(x^{N*},\xi^{N*})$ be the optimal solution and $LB_N$ be the optimal value.

\STATE \textbf{Feasibility-Check}

Solve the \textbf{FCP} \eqref{eq:feasibility-2} with $(x^{N*}, \xi^{N*})$. If $S_f(x^{N*},\xi^{N*})=0$, go the Step 4; otherwise, let $\hat w^{N*}$ be the worst-case scenario and let $UB_N=UB_{N-1}$, go to Step 5.

\STATE \textbf{Solve the Subproblem}

Solve the $\textbf{SP}$ \eqref{eq:subproblem} with $(x^{N*},\xi^{N*})$. Let $(w^{N*},\hat w^{N*}, z^{N*}, y^{N*})$ be the optimal solution and $UB_N=c^\top x^{N*}+g(\xi^{N*})+d^\top y^{N*}$, go to Step 5.

\STATE \textbf{Retrieve the active constraints}

Reduce redundant constraints (see Appendix \ref{apen-1}).

Deal with degenerate cases (see Section \ref{sec-3-2}).

Identify the set of active constraints $\mathcal{AC}_{N}$ based on \eqref{eq:active}.

\STATE if $|UB_N-LB_N|\le \epsilon$ and $S_f(x^{N*},\xi^{N*})=0$, terminate and output $(x^{N*},\xi^{N*})$. Otherwise, $N=N+1$, go to Step 2.

\end{algorithmic}
\end{algorithm} 

\begin{theorem}
\label{thm-1}
The Algorithm 1 generates the optimal solution to problem \eqref{eq:robust-whole} in finite iterations.
\end{theorem}

The proof of Theorem \ref{thm-1} can be found in Appendix \ref{apen-2}. The necessity and novelty of the proposed algorithm are highlighted below by comparing with two related algorithms.

\begin{itemize}
    \item[1)] \textbf{Nested C\&CG} \cite{zhao2012exact}. If we let $\hat w$ be the variable of the middle ``max'' problem and move the constraint $w=\min\{\hat w,\xi\}$ to the inner ``min'' problem. The problem obtained is equivalent to the original robust optimization model \eqref{eq:robust-whole} while the new uncertainty set becomes decision-independent. The nonlinear constraint $w=\min\{\hat w,\xi\}$ in the inner ``min'' can be linearized as
\begin{align}
   & w= \hat w z + \xi (1-z) \nonumber\\
 &    w \le \hat w,~ w \le \xi, z \in \{0,1\}^{N_w\times T}
\end{align}

Therefore, the problem is turned into a robust optimization problem with DIU and integer variables in the second stage. Denote the optimal value of the inner ``min'' as $Q(x,\hat w)$. Nested C\&CG \cite{zhao2012exact} might be applied to solve this problem.  However, since $Q(x,\hat w)$ is not quasiconvex in $\hat w$, there is no guarantee that the optimal solution of the ``max-min'' will reside at one extreme point (vertex) of the DIU set. Therefore, the algorithm may not terminate in finite steps. In case studies, we use nested C\&CG for comparison. The results show that, apart from the convergence issue, our method can greatly reduce the computational time compared with nested C\&CG. 

\item[2)] \textbf{Modified Bender's decomposition} \cite{zhang2021robust}. The method in \cite{zhang2021robust} deals with the robust optimization problem with a linear decision-dependent uncertainty set. Instead of adding the worst-case scenario to the master problem in each iteration, the corresponding dual variables are added with constraints to ensure that $w$ will be the optimal solution of the inner ``max-min'' problem with the selected dual variable. This constraint is equivalently represented as the KKT condition of the middle ``max'' problem so that the master problem renders a mixed integer linear program. However, for the model in this paper, the uncertainty set is a mixed-integer linear set, and we cannot get its KKT condition. 
\end{itemize}

\textbf{Remark}: The proposed algorithm can also be applied to solve a robust optimization under DDU with integer variables in the inner level. The integer variables in the inner level will only affect the solution methodology for the ``max-min'' \eqref{eq:max-min} (and the corresponding feasibility-check problem). Denote the integer variables by $\tau$, whose feasible set is $\Upsilon$. Then the ``max-min'' problem can be formulated as
\begin{align} \label{eq:stproblem}
    \mathop{\max}_{w \in \mathcal{W}(\xi)} ~ & \min_{\tau, y} d^{\top} y \nonumber\\
    \mbox{s.t.}~ & Fy \le f-Rw-Dx-\tilde E \tau, ~\tau \in \Upsilon
\end{align}

Instead of replacing the inner minimization problem with its KKT condition, we can separate the integer variable from the continuous variable and turn the problem \eqref{eq:stproblem} to the following tri-level optimization problem.
\begin{align} \label{eq:stproblem2}
    \mathop{\max}_{w \in \mathcal{W}(\xi)} ~ & \min_{\tau \in \Upsilon} ~ \min_y d^{\top} y \nonumber\\
    \mbox{s.t.}~ & Fy \le f-Rw-Dx-\tilde E \tau
\end{align}

For the problem \eqref{eq:stproblem2}, take $w$ as the first-stage decision variable, $\tau$ as the uncertain factor and $\Upsilon$ as the uncertainty set, and $y$ as the second-stage decision variable, then the C\&CG algorithm can be embedded to solve it.

\section{Illustrative Experiments}
\label{sec-4}
We first use the IEEE-39 bus system with a single period for demonstration. Then the IEEE-118 bus system and the real Guangdong power grid of China are tested to show the scalability and practicability of the proposed method. Detailed data can be found in \cite{Data}.

\subsection{IEEE-39 bus system}
\begin{figure}[t]
\centering
\includegraphics[width=0.9\columnwidth]{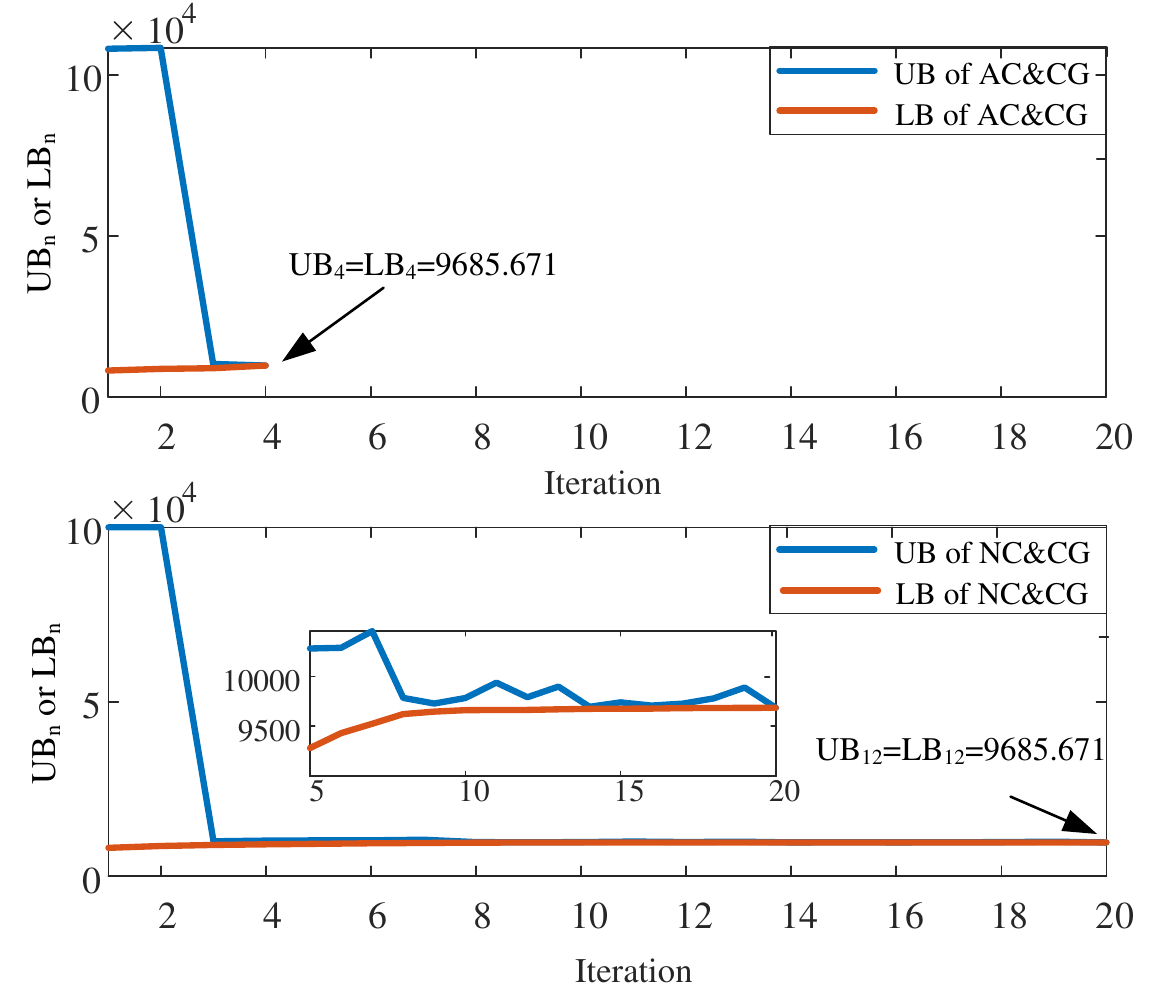}
\caption{Values of $UB_N$ and $LB_N$ during iterations.}
\label{fig:LBUB}
\end{figure}

The IEEE-39 bus system has 10 thermal generators and 46 transmission lines. The parameters of the thermal generators are given in TABLE \ref{tab:para}. There are three RGs connected to the power grid at nodes 4, 14, 29, respectively. Their forecast outputs are $W^e=[150, 150, 150]^\top$ MW and $W^l=0.5W^e, W^u=1.5W^e$. The uncertainty budget $\Gamma_S=2$. The proposed Adaptive C\&CG (AC\&CG) algorithm and the Nested C\&CG (NC\&CG) algorithm are applied to solve the RGD problem, respectively. The change of $UB_N$ and $LB_N$ during iterations are given in Fig. \ref{fig:LBUB}. Note that the number of iterations of NC\&CG in Fig. \ref{fig:LBUB} and TABLE \ref{tab:uncertain} refer to the number of outer-loops for adding the worst-case scenarios to the master problem. We can find that both algorithms output the same solution with an optimal objective value of \$9,685.671 and the curtailment strategy $\xi=[195.00, 200.00, 175.00]^\top$ MW. This verifies the proposed algorithm. The NC\&CG algorithm takes 41.79s to converge while the AC\&CG takes 2.61s. The proposed AC\&CG algorithm can greatly reduce the computational time since it only requires 4 iterations to reach the optimal solution whereas the NC\&CG takes 20 iterations. Let alone the fact that in each iteration, the NC\&CG algorithm takes another 3-4 iterations to solve the subproblem.
\begin{table}[t]
        \renewcommand{\arraystretch}{1.3}
        \renewcommand{\tabcolsep}{1em}
        \scriptsize
        \centering
        \caption{Parameters of the thermal generators}
        \label{tab:para}
        \begin{tabular}{ccccccc}
                \hline 
               Unit & $P_i^l$ & $P_i^u$ & $R_i^{\pm}$ & $\alpha_i$ & $\beta_i^{\pm}$ & $d_i^{\pm}$\\
               & (MW) & (MW) & (MW) & (\$/MW) & (\$/MW) & (\$/MW)\\
               \hline
               1 & 70 & 155 & 120 & 16.19 & 1.62 & 1.62 \\
               2 & 70 & 155 & 120 & 17.26 & 1.73 & 1.73 \\
               3 & 40 & 80 & 60 & 16.60 & 1.66 & 1.66\\
               4 & 30 & 50 & 50 & 16.50 & 1.59 & 1.59 \\
               5 & 55 & 100 & 70 & 19.70 & 1.97 & 1.97 \\
               6 & 30 & 50 & 30 & 22.26 & 2.22 & 2.22\\
               7 & 25 & 55 & 34 & 22.74 &  2.77 & 2.77\\
               8 & 20 & 55 & 22 &  25.92 & 2.59 & 2.59\\ 
               9 & 20 & 55 & 20 & 27.27 & 2.73 & 2.73 \\
               10 & 20 & 55 & 10 & 27.79 & 1.30 & 1.30\\
               \hline
        \end{tabular}
\end{table}

\begin{figure*}[t]
\centering
\includegraphics[width=1.8\columnwidth]{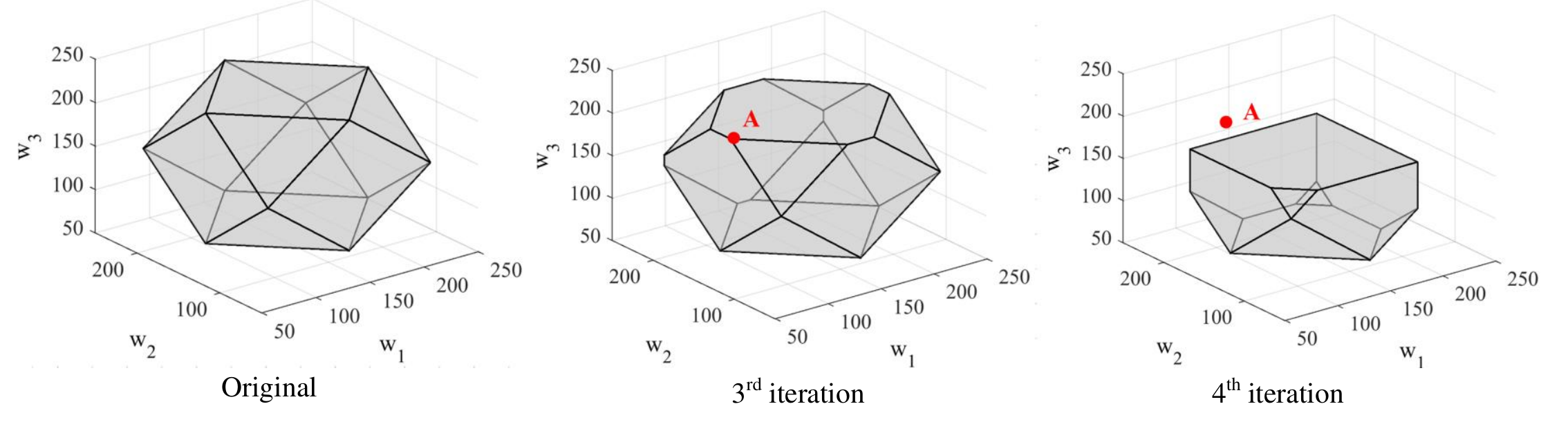}
\caption{Decision-dependent uncertainty sets.}
\label{fig:DDU}
\end{figure*}

\begin{figure}[t]
\centering
\includegraphics[width=0.8\columnwidth]{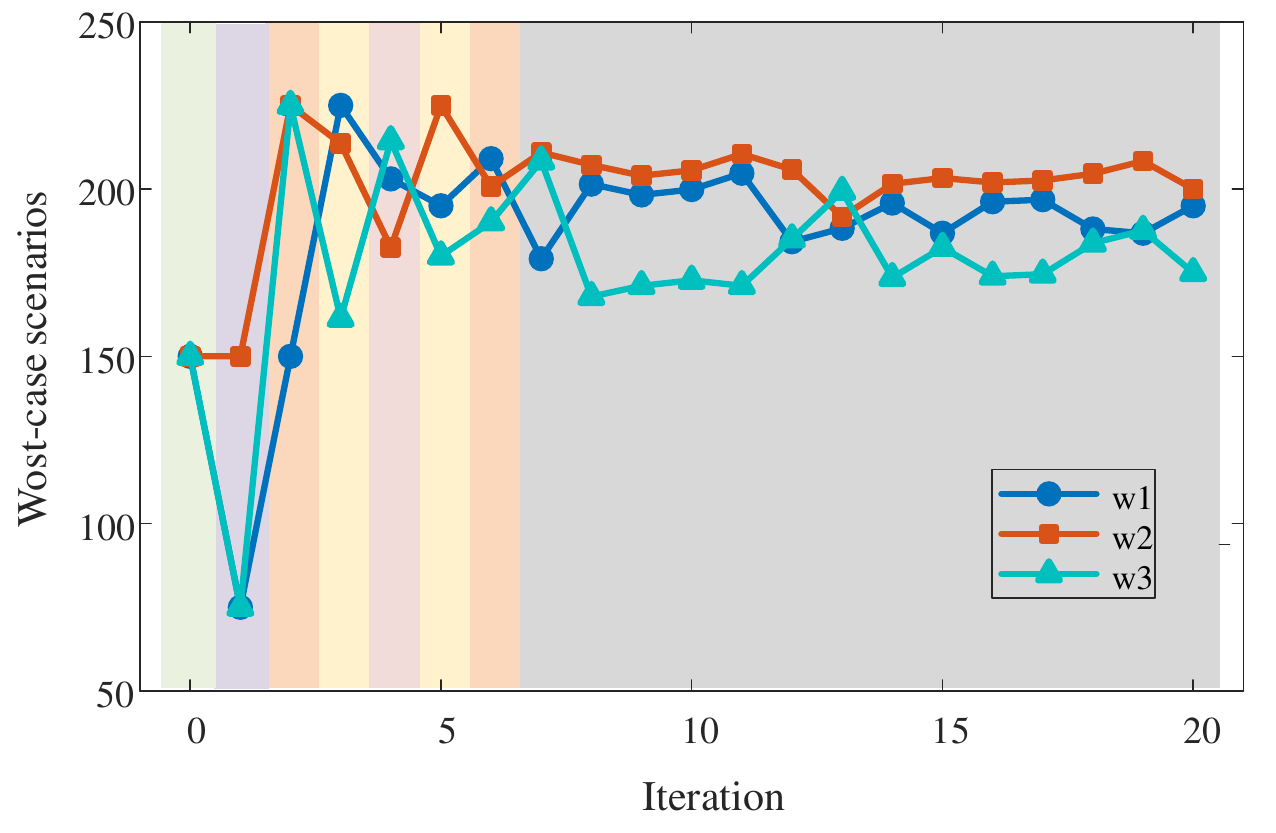}
\caption{Worst-case scenarios during iterations of the nested C\&CG.}
\label{fig:NCCG}
\end{figure}

\begin{figure}[t]
\centering
\includegraphics[width=0.8\columnwidth]{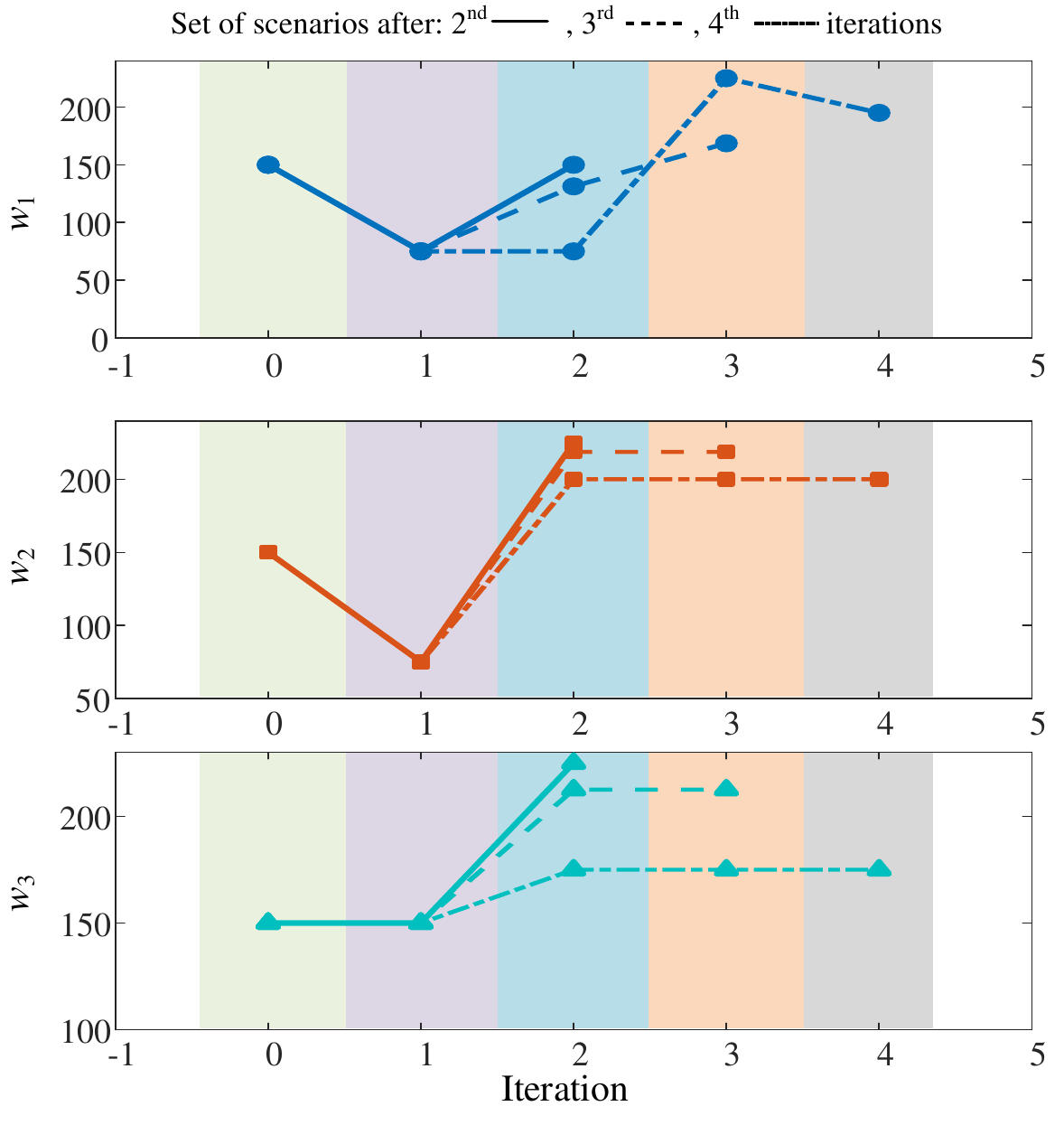}
\caption{Worst-case scenarios during iterations of our proposed algorithm.}
\label{fig:MCCG}
\end{figure}

The decision-dependent uncertainty set in the 1st, 3rd, and 4th iterations are shown in Fig. \ref{fig:DDU}. When the first-stage decision changes, the uncertainty set changes dramatically. As a result, a previously picked scenario, for example if it is point A, may fall outside of the uncertainty set. Continuing to use traditional robust optimization algorithms may result in over-conservative or infeasible results. This demonstrates the necessity of this work. The worst-case scenarios selected during the iterations of the two algorithms are plotted in Fig. \ref{fig:NCCG} and Fig. \ref{fig:MCCG}, respectively. In each iteration, a different background color represents a different set of active constraints.

In Fig. \ref{fig:NCCG}, the same set of active constraint recurs. For example, in the 7th - 20th iterations, the worst-case scenarios are generated by the same active constraint $w=\xi$. This is the root of the inefficiency of traditional algorithms: when the first-stage decision changes, the previously selected scenario may no longer be a vertex of the new uncertainty set. Therefore, there are an infinite number of possible uncertain scenarios and the algorithm does not necessarily stop in a finite number of steps. In contrast, in Fig. \ref{fig:MCCG}, the same set of active constraints appears only once. Instead of a fixed value, the previously selected scenario will change with the first-stage decision. The set of selected scenarios after the 2nd, 3rd, and 4th iterations are connected by the solid line, dashed line, and dash-dot line, respectively. For example, the scenario identified in the 2nd iteration is $w^{2*}=[150, 225, 225]^\top$ MW, and thus, the active constraints are $w^2_{2}=\xi_2$, $w_{3}^2=\xi_3$, and $w_1^2=w_2^{2}+w_2^{3}-300$. In the 3rd iteration, it changes to $w_2^*=[131.34, 218.78, 212.56]^\top$ since the curtailment strategy changes to $\xi=[225, 218.78, 212.56]^\top$ MW. Then as the curtailment strategy changes to $\xi=[195.00, 200.00, 175.00]^\top$ MW in the 4th iteration, the scenario further changes to $w_2^*=[75.00, 200.00, 175.00]$ MW.

To show the advantages of the proposed model (Model-1), we compare it with the one that merely uses re-dispatch curtailment (Model-2). The Model-2 is derived by fixing $\xi_{jt}$ in the proposed model to $W_{jt}^u$ for all $j\in [N_w], t \in [T]$. Since $\xi$ is now a constant, Model-2 is a robust optimization with decision-independent uncertainty. We test the performance of the two models under different maximum reserve limits. The results are shown in TABLE \ref{tab:comparison-1}. We can find that with the proposed model, the reserve cost can always be reduced. For example, in the case with $R_i^{\pm}$, the reserve cost is decreased by $(461.9-442.2)/461.9=4.3\%$. Moreover, the real-time curtailment in the proposed model is much smaller than that in Model-2, so the system dynamic characteristics do not change significantly, reducing the risk of instability. The total cost of the proposed model is substantially lower than that of Model-2, and with a more stringent reserve limit, the cost reduction becomes more significant.

\begin{table*}[t]
        \renewcommand{\arraystretch}{1.3}
        \renewcommand{\tabcolsep}{1em}
        \scriptsize
        \centering
        \caption{Comparison of the proposed model and Model-2}
        \label{tab:comparison-1}
        \begin{tabular}{ccccccc}
                \hline 
                & \multicolumn{2}{c}{$R_i^{\pm}$} & \multicolumn{2}{c}{$0.5R_i^{\pm}$} & \multicolumn{2}{c}{$0.3R_i^{\pm}$} \\
              & Model-1 & Model-2 & Model-1 & Model-2 & Model-1 & Model-2\\
              \hline
              Generation cost (\$) & 8,129.3 & 8,129.3 & 8,143.0 & 8,149.6 & 8,215.7 & 8,215.7\\
              Reserve cost (\$)	& 442.2 & 461.9 & 450.8 & 468.1 & 465.6 & 465.6\\
              Pre-dispatch curtailment cost (\$) &	787.5 & 0 & 787.5 & 0 & 787.5 & 0\\
              Redispatch cost (\$) & 183.9 & 177.9 & 179.3 & 180.0 & 187.5 & 187.5\\
              Re-dispatch curtailment Amount (MW) / Cost (\$) & 4.76/142.7 & 39.91/1,146.0 &	9.91/148.7 & 39.91/1,146.0 & 9.91/148.7	& 39.91/1,146.0\\
              Total cost (\$) & 9,685.6 & 9,915.1 & 9,709.3 & 9,943.7 & 9,805.0 & 10,014.8\\
               \hline
        \end{tabular}
\end{table*}

Out-of-sample test is also conducted to evaluate the performance of the proposed model in terms of feasibility and cost. We randomly generate scenarios from a normal distribution with the mean of $W_{jt}^e$ and the standard variance of 30, 40, 50, respectively. Five hundred scenarios are tested for each setting and the number of infeasible scenarios, the expected re-dispatch cost, and the worst-case real-time curtailment are recorded in TABLE \ref{tab:comparison-2}. We can find that the proposed model can achieve a similar performance to the Model-2 in robustness but with a lower expected re-dispatch cost. Moreover, the worst-case real-time curtailment of the proposed model is much smaller than that of Model-2, meaning that the proposed model is posing less pressure on the real-time power system operation.

\begin{table*}[t]
        \renewcommand{\arraystretch}{1.3}
        \renewcommand{\tabcolsep}{1em}
        \scriptsize
        \centering
        \caption{Out-of-sample test of the proposed model and Model-2}
        \label{tab:comparison-2}
        \begin{tabular}{ccccccc}
                \hline 
               Standard variance & \multicolumn{2}{c}{30} & \multicolumn{2}{c}{40} & \multicolumn{2}{c}{50} \\
              & Model-1 & Model-2 & Model-1 & Model-2 & Model-1 & Model-2\\
              \hline
              Infeasible scenarios & 22/500 & 25/500 & 42/500 & 46/500 & 55/500 & 68/500\\
              Expected re-dispatch cost (\$) & 106.74 & 107.12 & 111.10 & 116.28 & 120.17 & 127.46\\
              Worst-case real-time curtailment (MW) & 9.91 & 173.59 & 9.91 & 308.30 & 9.91 & 371.61\\
               \hline
        \end{tabular}
\end{table*}

\subsection{IEEE-118 bus system: Benchmark}
\begin{figure}[t]
\centering
\includegraphics[width=1.0\columnwidth]{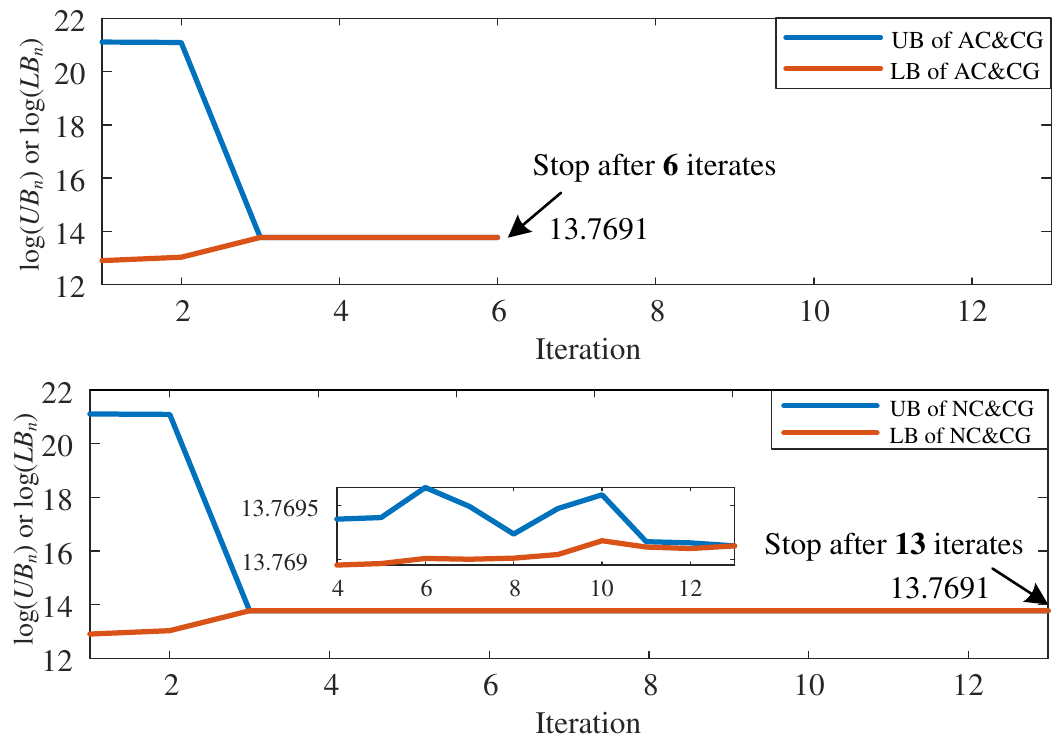}
\caption{Values of \textbf{log}($UB_N$) and \textbf{log}($LB_N$) during iterations.}
\label{fig:UBLB118}
\end{figure}

We further test the performance of our proposed algorithm using a larger IEEE-118 bus system. There are 54 thermal generators and 186 transmission lines. Our proposed algorithm (AC\&CG) takes 6 iterations (35.17 seconds) to reach the robust optimal solution. The traditional nested C\&CG takes 13 iterations (418.22 seconds). The change of $UB_N$ and $LB_N$ are given in Fig. \ref{fig:UBLB118}. Since $UB_N$ and $LB_N$ vary in a wide range, therefore, we draw the \textbf{log}($UB_N$) and \textbf{log}($LB_N$) instead. Both algorithms output the same result while the proposed algorithm can greatly reduce the computational time.

\begin{figure}[t]
\centering
\includegraphics[width=0.8\columnwidth]{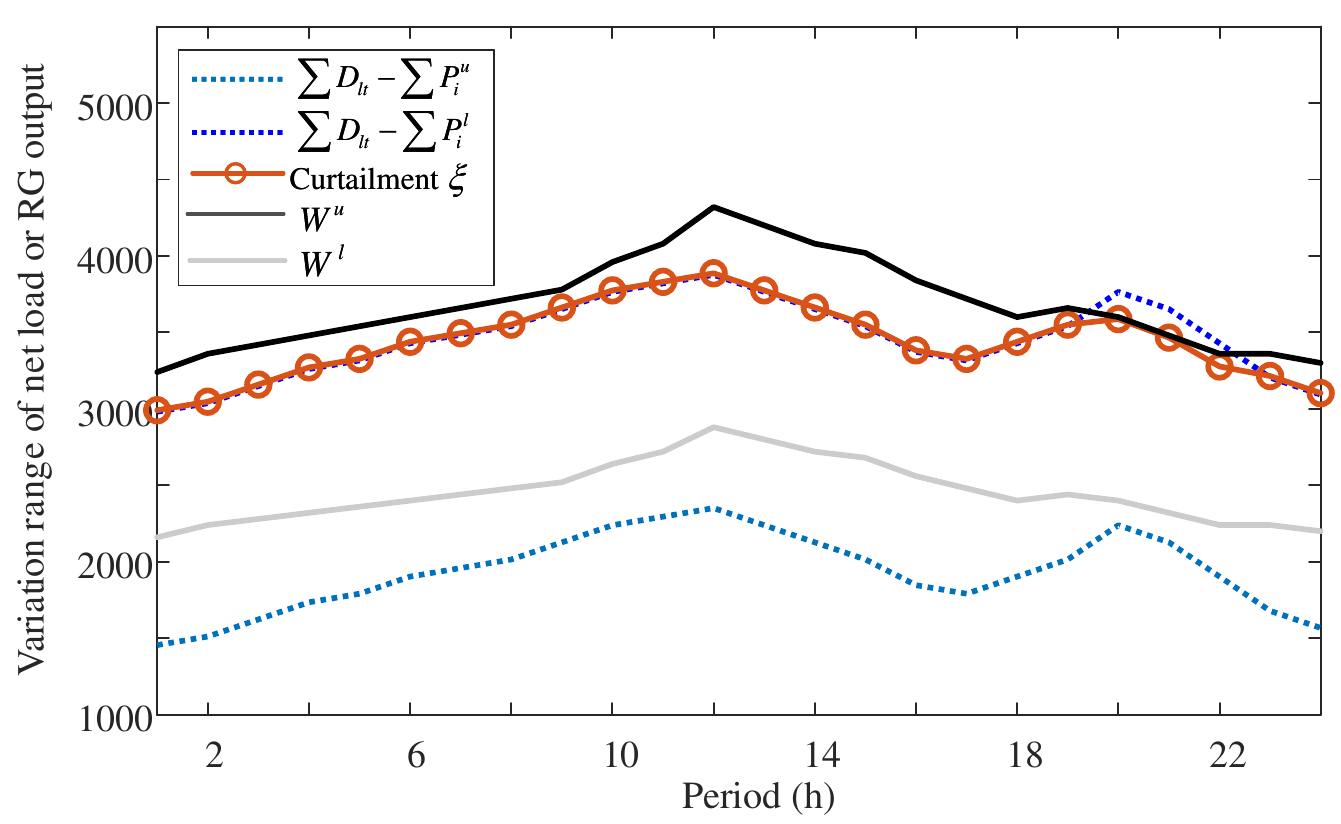}
\caption{Curtailment schedule and variation range of net load and RG output.}
\label{fig:Result118}
\end{figure}

The optimal curtailment strategy ($\xi^*$) as well as the upper/lower bound of the net load ($\sum_l D_{lt} - \sum_i p_{it}$) and the RG power output ($\sum_j w_{jt}$) is shown in Fig. \ref{fig:Result118}. With growing renewable energy penetration, in many periods (e.g., 1-19h), the upper bound of the total RG power output is much higher than the upper bound of the net load, i.e., $W^u > \sum_l D_{lt}-\sum_i P_i^{l}$. Therefore, preparatory curtailment is required to avoid potential risks brought by a large amount of real-time curtailment; and correspondingly, the red curve in Fig. \ref{fig:Result118} is lower than the upper black curve ($W^u$). In hours 20 and 21, the maximum net load is higher than the maximum RG power output, and there is no need of preparatory curtailment. However, in hour 22, even though the maximum net load is higher than the maximum RG power output, pre-dispatch curtailment still occurs, which is caused by the ramping capacity limit. From Fig. \ref{fig:Result118} we can find that, curtailment happens in most of the time, and if all are done in real-time without a preparatory schedule, it will result in frequent emergency controls that will jeopardize the system's stability. That's why we need to decide on the curtailment schedule in the pre-dispatch stage.

\begin{figure}[t]
\centering
\includegraphics[width=1.0\columnwidth]{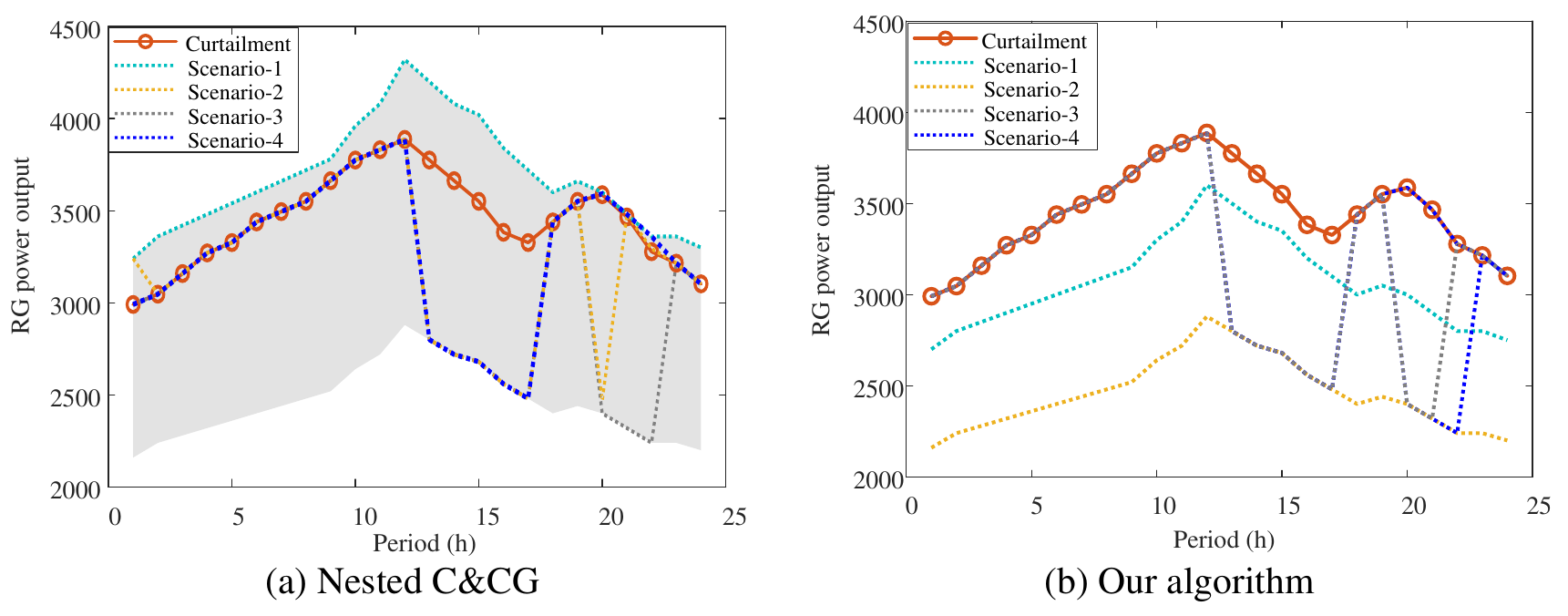}
\caption{Selected worst-case scenarios under two algorithms.}
\label{fig:worstcase}
\end{figure}

Four typical scenarios are selected from the uncertainty set in the last iteration of the nested C\&CG algorithm and the proposed algorithm, respectively. They are plotted in Fig. \ref{fig:worstcase} together with the optimal pre-dispatch curtailment strategy $\xi$. When applying the nested C\&CG algorithm, we move the $w=\min\{\hat w, \xi\}$ constraint to the inner ``min'' problem and turn the uncertainty set into a decision-independent one, which is the shaded area in Fig. \ref{fig:worstcase}(a) (for more details, please refer to Section \ref{sec-3-4}). First, we can find that scenarios 1, 2, and 4 exceed the upper limit imposed by the curtailment schedule. Thus, they are not feasible for the original decision-dependent uncertainty set \eqref{eq:DDU-set} and this makes the traditional C\&CG not applicable. Moreover, the selected scenarios (e.g. scenarios 2 and 3) may not be vertices of the decision-independent uncertainty set (shaded area). Therefore, there is no guarantee that the nested C\&CG algorithm will stop after a finite number of iterations. In contrast, under the proposed algorithm, the previously selected algorithm can change adaptively and all remains feasible for the new uncertainty set as in Fig. \ref{fig:worstcase}(b).

\subsection{IEEE-118 bus system: Impact of different factors}
In the following, we analyze the impact of several factors. First, we change the ramping limit $\mathcal{R}_i^{u}, \mathcal{R}_i^d, \forall i$ from 1.0 to 3.0 times of their original values. The results are shown in TABLE \ref{tab:ramp}. With a tighter limit, the system operator has to use more costly generators for reserve, so the reserve cost increases. To cope with this, we want to reduce the uncertainty faced by the system, so that less reserve will be needed. The level of uncertainty can be reduced by planning most of the curtailment in advance. For example, the total curtailment in the case with ratio=1.0 is similar to that in the case with ratio=3.0. But in the former case, more curtailment are prepared beforehand. With more renewable generators replace the conventional generators, the ramping capacity will become more scarce, and thus, preparatory curtailment is more imperative.
When the ratio changes from 2.0 to 2.5, even though the curtailments are the same, the total cost declines. This is because larger ramping limits enable a better allocation of contemporary output and reserve capacity of thermal units. When the ratio changes from 2.5 to 3.0, both the curtailment cost and total cost remain unchanged since ramping limits are no longer binding factors. Moreover, the computational times are all acceptable.

\begin{table}[t]
        \renewcommand{\arraystretch}{1.3}
        \renewcommand{\tabcolsep}{1em}
        \centering
        \footnotesize
        \caption{Costs and times under different ramping limits}
        \label{tab:ramp}
        \begin{tabular}{cccccc}
                \hline 
              Ratio & First/Second stage   & reserve cost & total cost & Time\\
              & curtailment & & & \\
              \hline
              1.0 & 5,610.4 / 160.0 & 32,668.5 & 954,675.6 & 35.17\\
              1.3 &  5,543.3 / 170.0 & 31,195.9 & 936,187.3 & 46.74\\
                1.5 & 5,533.3 / 214.0 & 30,483.7 & 930,556.1 & 39.60\\
                2.0 & 5,533.9 / 210.0 & 28,626.8 & 924,994.1 & 18.57\\
                2.5 &  5,533.9 / 210.0 & 27,961.5 & 924,266.8 & 22.04\\
                3.0 & 5,533.9 / 210.0 & 27,961.5 & 924,266.8 & 29.47\\
                \hline
        \end{tabular}
\end{table}


Furthermore, we test the performance of the two algorithms under different levels of uncertainty as given in TABLE \ref{tab:uncertain}. When the uncertainty decreases, both the curtailment amount and the total cost decline, which is straightforward. With higher uncertainty, the number of iterations needed by the proposed AC\&CG algorithm is stable while that taken by the NC\&CG is less predictable. This shows the superiority of the proposed algorithm in the future power system with severe renewable energy fluctuations. 

\begin{table}[t]
        \renewcommand{\arraystretch}{1.3}
        \renewcommand{\tabcolsep}{0.5em}
        \centering
        \footnotesize
        \caption{Curtailments and times under different uncertainty levels}
        \label{tab:uncertain}
        \begin{tabular}{cccc}
                \hline 
               $[W^l, W^u]$ & First/Second stage & \multicolumn{2}{c}{AC\&CG / NC\&CG} \\
               & curtailment & No. of iterates & Time (s)\\
               \hline
            $[0.80, 1.20]W^e$ & 5,610.4 / 160.0 &  6 / 13 & 35.17 / 418.22\\
            $[0.85, 1.15]W^e$ & 2,646.3 / 93.2 &  5 / 27 & 46.19 / 861.09\\
            $[0.90, 1.10]W^e$ & 885.6 / 23.7 &  7 / 17 & 82.00 / 240.46\\
            $[0.95, 1.05]W^e$ & 142.7 / 84.3 &  3 / 3 & 16.53 / 23.18\\
            $[1.00, 1.00]W^e$ & 0.0 / 0.0 &  1 / 1 & 2.01 / 3.49\\
                \hline
        \end{tabular}
\end{table}

\subsection{Guangdong Power System}
\begin{figure}[t]
\centering
\includegraphics[width=0.75\columnwidth]{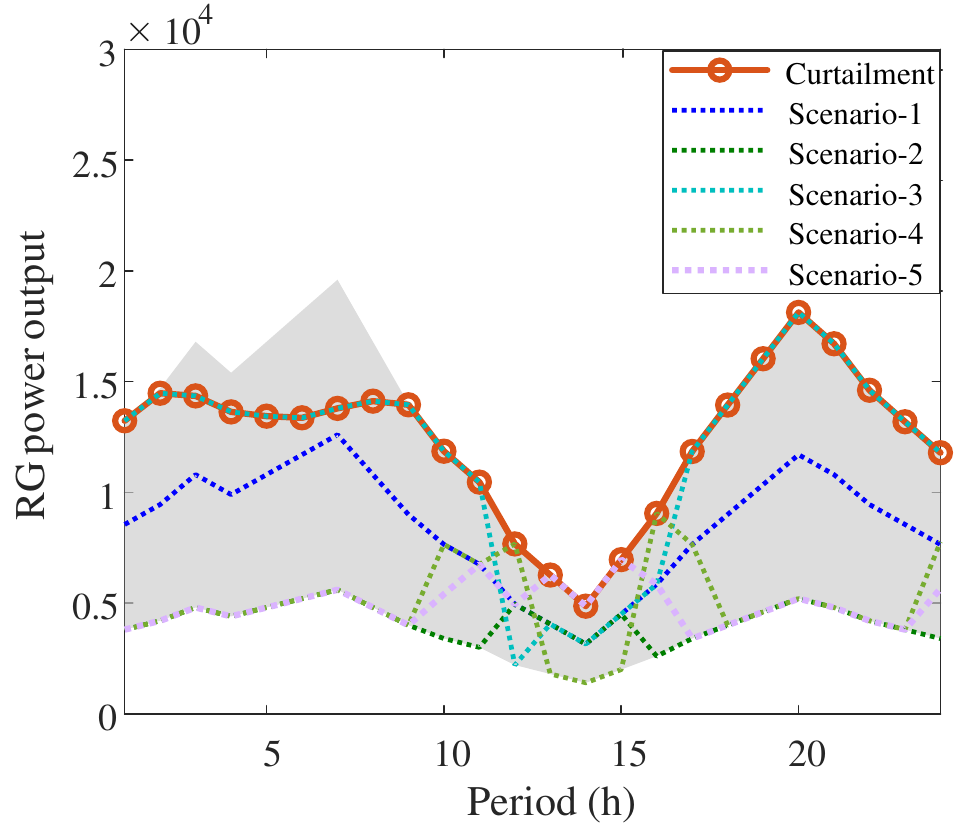}
\caption{Optimal pre-dispatch curtailment and worst-case scenarios of the Guangdong system.}
\label{fig:guangdong}
\end{figure}

The real Guangdong power grid of China is also tested to show the scalability of the proposed approach. The whole system has 1880 buses, 2452 transmission lines, 174 thermal generators, 453 loads, and six in planning large-scale wind farms located at Guangzhou, Shaoguan, Shenzhen, Dongguan, Shantou, and Zhanjiang, respectively. The predicted nodal loads in a normal winter day are used in this test. We change the uncertainty budget $\Gamma_T$ from 18 to 24 and record the total cost, number of iterations, and computational times in TABLE \ref{tab:guangdong}. We can find that the computational times are less than 18min, which is acceptable. The optimal pre-dispatch curtailment and the worst-case scenarios of the case with $\Gamma_T=20$ are shown in Fig. \ref{fig:guangdong}. With the pre-dispatch curtailment $\xi$, the variation range of renewable power output (capped by $\xi$) is smaller than the original region $[W^l, W^u]$ (the grey area). All selected scenarios remain within the decision-dependent uncertainty set $\mathcal{W}(\xi)$, which validates the theoretical analysis.

\begin{table}[t]
        \renewcommand{\arraystretch}{1.3}
        \renewcommand{\tabcolsep}{0.5em}
        \centering
        \scriptsize
        \caption{Total cost and computational time of the Guangdong system}
        \label{tab:guangdong}
        \begin{tabular}{cccccc}
        \hline
        $\Gamma_T$ & 16 & 18 & 20 & 22 & 24\\
                \hline 
           Total cost (\$) & 170,416,392 & 170,660,423 & 170,864,676	& 171,019,459 & 171,134,996\\
           Iteration & 6 & 6 & 5 & 5 & 3\\
           Time (s) & 1075.1 & 911.0 & 591.8 & 154.1 & 43.2\\
                \hline
        \end{tabular}
\end{table}

\section{Conclusion}
\label{sec-5}
In this paper, a novel RGD model considering preparatory schedule of renewable power curtailment is proposed. The model casts down to a robust optimization with DDU. An adaptive C\&CG algorithm was developed to solve the problem efficiently by identifying the active constraints in each iteration and adding them to the master problem. In this way, the worst-case scenarios remain to be vertices of the uncertainty set when the first-stage decision changes. Our main findings are:

1) We prove theoretically that the proposed algorithm can generate the robust optimal solution in a finite number of iterations. Many of the existing algorithms do not have this guarantee due to the decision-dependent uncertainty.

2) The proposed algorithm can reduce the computational time by 95\% compared with the nested C\&CG algorithm.

3) With higher uncertainty, the performance of the proposed algorithm remains stable while the number of iterations needed by the nested C\&CG is less predictable.

4) Compared with the traditional model with real-time curtailment only, the proposed model can avoid significant wind curtailments in the re-dispatch stage and can also lower the amount of reserve required.


\appendix
\makeatletter
\@addtoreset{equation}{section}
\@addtoreset{theorem}{section}
\makeatother
\setcounter{equation}{0}
\renewcommand{\theequation}{A.\arabic{equation}}
\subsection{Redundancy Elimination}
\label{apen-1}
A redundant constraint in $\mathcal{W}(\xi)$ in defined as follows.

\vspace{0.5em}
\noindent \textbf{Definition A1.} (Redundant Constraint)
Given $\xi$ and $z^*$, a constraint is \emph{redundant} if removing that constraint does not change the uncertainty set $\mathcal{W}(\xi)$. 
\vspace{0.5em}

Denote $Q:=[H,-C]$, $u:=[w^\top, \hat w^\top]^\top$, and $q:=F\xi+Gz^*+a$. Then the redundant constraints can be removed with the help of the following theorem.

\vspace{0.5em}
\noindent \textbf{Theorem A1.} Constraint $Q_j u \le q_j$ is redundant if the following linear program has a feasible solution
\begin{align}
    v^\top  Q_{[-j]} = Q_j,~ q_j \ge v^\top q_{[-j]},~ v \ge 0
\end{align}

\begin{proof}
If constraint $Q_j u \le q_j$ is redundant, the polytope $\{u|Q_{[-j]} u \le q_{[-j]}\}$ defined by the remaining constraints must be a subset of $\{u|Q_ju \le q_j\}$. Therefore, we have $\{u|Q_{[-j]}u \le q_{[-j]}\} \cap \{u|Q_ju>q_j\}=\emptyset$. According to the Nonhomogeneous Farkas Lemma \cite{xu1994generalization}, it is equivalent to 
\begin{align}
    \mathcal{Q}=\{v|v^\top Q_{[-j]}=Q_j,~ q_j \ge v^\top q_{[-j]}, ~ v \ge 0\} \ne \emptyset
\end{align}
In other words, the corresponding linear program has a feasible solution.
\end{proof}

By checking each of the constraints in $\mathcal{W}(\xi)$, remove the redundant ones, and identify the active constraints, we can avoid the primal degeneracy.

\setcounter{equation}{0}
\renewcommand{\theequation}{B.\arabic{equation}}
\subsection{Proof of Theorem \ref{thm-1}}
\label{apen-2}
We start the proof of Theorem \ref{thm-1} with the following lemma.

\vspace{0.5em}
\noindent \textbf{Lemma B1}. Suppose the optimal solution of \eqref{eq:robust-whole} is $(x^*,\xi^*)$. Let $N$ be the iteration run of Algorithm 1. For any $n \in [N]$:

\begin{itemize}
    \item[(a)] $LB_n \le c^\top x^*+g(\xi^*) +S(x^*,\xi^*) \le UB_n$
    \item[(b)] For any $n_1, n_2 \in [N-1]$, the set of active constraints $\mathcal{AC}_{n_1}$ won't be the same as the set of active constraints $\mathcal{AC}_{n_2}$.
    \item[(c)] If there exists some $n \in [N-1]$ whose set of active constraints $\mathcal{AC}_{n}$ is the same as $\mathcal{AC}_N$, the algorithm terminates in $N$.
\end{itemize}

\begin{proof} 

\textbf{Assertion (a):} For each $n \in [N]$, constraints \eqref{eq:master-4} and \eqref{eq:master-5} give a feasible point of $\mathcal{W}(\xi)$. The robust optimization problem \eqref{eq:robust-whole} can be equivalently written as
\bsq
\label{eq:robusteq}
\begin{align}
    \min_{x,\xi}~ & c^\top x+g(\xi)+\eta \label{eq:robusteq-1}\\
    \mbox{s.t.}~ & x \in \mathbb{X}, \xi \in \Xi \label{eq:robusteq-2}\\
    ~ & \eta \ge Q(x,w),~\forall w \in \mathcal{W}(\xi) \label{eq:robusteq-3}
\end{align}
\esq
while the master problem \textbf{MP} is
\bsq
\label{eq:mastereq}
\begin{align}
    \min_{x,\xi}~ & c^\top x+g(\xi)+\eta \label{eq:mastereq-1}\\
    \mbox{s.t.}~ & x \in \mathbb{X}, \xi \in \Xi \label{eq:mastereq-2}\\
    ~ & \eta \ge Q(x,w),~\forall w \in \Gamma(\textbf{Sol}(\mathcal{AC}_n,\xi)) \label{eq:mastereq-3}
\end{align}
\esq
where $\textbf{Sol}(\mathcal{AC}_n,\xi)$ refers to the unique solution of the set of active constraints $\mathcal{AC}_n$ given $\xi$ and $\Gamma(.)$ refers to the projection exerted by \eqref{eq:master-5}. Then, $w^n$ is a feasible point of $\mathcal{W}(\xi)$. The two problems \eqref{eq:robusteq} and \eqref{eq:mastereq} have the same objective functions while the constraint \eqref{eq:mastereq-3} is a relaxation of \eqref{eq:robusteq-3}. Therefore, the optimal value of \eqref{eq:mastereq}, i.e., $LB_n$, is no more than that of \eqref{eq:robusteq}, i.e., $c^\top x^*+g(\xi^*)+S(x^*,\xi^*)$. Next, we prove $UB_n \ge c^{\top}x^*+g(\xi^*)+S(x^*,\xi^*)$ by induction. First of all, $UB_0>c^{\top}x^*+g(\xi^*)+S(x^*,\xi^*)$. Suppose for the sake of induction that $UB_{n-1}>c^{\top}x^*+g(\xi^*)+S(x^*,\xi^*)$. Then if $S_f(x^{n*},\xi^{n*})>0$, we have $UB_{n}=UB_{n-1}>c^{\top}x^*+g(\xi^*)+S(x^*,\xi^*)$; else, $(x^{n*},\xi^{n*})$ is robust feasible and  we have $UB_n=c^\top x^{n*}+g(\xi^{n*})+S(x^{n*},\xi^{n*})$. Due to the optimality of $(x^*,\xi^*)$, we have $UB_n \ge c^\top x^*+g(\xi^*)+S(x^*,\xi^*)$.


\textbf{Assertion (b):} If we have $\mathcal{AC}_{n_1}$ be the same as $\mathcal{AC}_{n_2}$, then, $w^{n_2*}=\Gamma(\textbf{Sol}(\mathcal{AC}_{n_2},\xi^{n_2*}))=\Gamma(\textbf{Sol}(\mathcal{AC}_{n_1},\xi^{n_2*}))$. This also implies that $(x^{n_2*},\xi^{n_2*})$ is robust feasible, so we have $UB_{n_2}=c^\top x^{n_2*} + g(\xi^{n_2*}) + Q(x^{n_2*},w^{n_2*})$. Moreover, $\eta$ in the master problem is the maximum of $Q(x^{n_2*}, w)$ for all $w=\Gamma(\textbf{Sol}(\mathcal{AC}_j,\xi^{n_2*})), j=[n_2-1]$. so $\eta \ge Q(x^{n_2*},w^{n_1*})=Q(x^{n_2*},w^{n_2*})$. Therefore, $LB_{n_2} \ge UB_{n_2}$. Together with $UB_{n_2} \ge LB_{n_2}$, we have $LB_{n_2}=UB_{n_2}$, and the algorithm terminates in the $n_2$ iteration, which is contradict to the fact that $N>n_2$ is the iteration run.

\textbf{Assertion (c):} Following a similar procedure as in the proof of (b), it is easy to get (c).
\end{proof}

Now the proof of \textbf{Theorem \ref{thm-1}} is given below.

First, we prove that the algorithm terminates in a finite number of iterations. According to Lemma B1(b)-(c), the same set of active constraints won't be attached twice. Moreover, each set of active constraints is composed of $4N_w\times T$ linearly-independent constraints from $\mathcal{W}(\xi)$. Since there are a finite number of those possible combination of constraints, when A1 holds, the algorithm will terminate in finite steps.

Next, we show the robust feasibility of $(x^{N*},\xi^{N*})$. Since $(x^{N*},\xi^{N*})$ is generated by the \textbf{MP}, we have $x^{N*} \in \mathbb{X}$ and $\xi^{N*} \in \Xi$. Also, the algorithm terminates when $S_f(x^{N*},\xi^{N*})=0$, so for any $w \in \mathcal{W}(\xi^{N*})$, the second-stage problem \eqref{eq:SD} is feasible.

Finally, we show the optimality of $(x^{N*},\xi^{N*})$. According to Lemma B1(a), we have $LB_N \le c^\top x^*+g(\xi^*)+S(x^*,\xi^*)\le UB_N$. Together with the condition for termination $|UB_N-LB_N|\le \epsilon$. Therefore, 
\begin{align}
   ~ &  |c^\top x^N+g(\xi^N)+S(x^N,\xi^N)- (c^\top x^*+g(\xi^*)+S(x^*,\xi^*))| \nonumber\\
   = ~ & |UB_N -(c^\top x^*+g(\xi^*)+S(x^*,\xi^*))| \nonumber\\
   \le ~ & |UB_N-LB_N| \le \epsilon
\end{align}
which justifies the optimality of $(x^{N*},\xi^{N*})$.  $\hfill \square$

\ifCLASSOPTIONcaptionsoff
\newpage
\fi

\bibliographystyle{IEEEtran}
\bibliography{IEEEabrv,mybib}

\end{document}